\documentclass[reqno]{amsart}
\usepackage{nicematrix}
\usepackage{tikz}
\usepackage{xcolor}
\usepackage{amsmath,bm,bbm,amsthm, amssymb}
\usepackage{mathtools}
\usepackage{fullpage}
\usepackage{color}
\usepackage{dsfont}
\usepackage{animate}
\usepackage{etoolbox}
\usepackage{comment}
\usepackage{pgf}
\usetikzlibrary{calc}
\usetikzlibrary{patterns}
\usetikzlibrary{arrows}
\usetikzlibrary{decorations.pathreplacing}
\usepackage[utf8]{inputenc}
\usepackage{pgfplots}

\usepackage{tcolorbox}
\usepackage{adjustbox}
\usepackage[final]{pdfpages}
\usepackage[ruled,vlined,linesnumbered,algosection,resetcount]{algorithm2e}
\usepackage{blkarray}
\usepackage{arydshln}
\usepackage{wrapfig}
\usepackage[colorlinks=true]{hyperref}

\theoremstyle{plain}
\newtheorem{theorem}{Theorem}
\newtheorem{proposition}[theorem]{Proposition}
\newtheorem{corollary}[theorem]{Corollary}
\newtheorem{lemma}[theorem]{Lemma}

\theoremstyle{definition}
\newtheorem{definition}[theorem]{Definition}

\theoremstyle{remark}

\def\R{\mathbb R}

\def\Z{\mathbb Z}

\def\P{\mathbb P}

\def\E{\mathbb E}
\def\T{W}

\def\xx{\mathbf x}
\def\yy{\mathbf y}

\def\00{\mathbf 0}

\def\bet{\begin{theorem}}
\def\ent{\end{theorem}}
\def\bec{\begin{corollary}}
\def\enc{\end{corollary}}
\def\bep{\begin{proof}}
\def\enp{\end{proof}}
\def\bede{\begin{definition}}
\def\ende{\end{definition}}
\def\f{\frac}

\def\g{\gamma}
\def\la{\lambda}

\def\ms{\mathsf}
\def\co{\colon}
\def\N{\mathbb N}
\def\mc{ \mathcal}

\def\PP{\mc P}

\def\one{\mathds1}

\def\d{\mathrm d}

\renewcommand\leq{\leqslant}

\renewcommand\le{\leqslant}
\renewcommand\ge{\geqslant}

\def\bel{\begin{lemma}}
\def\enl{\end{lemma}}
\def\im{\item}
\def\been{\begin{enumerate}}
\def\enen{\end{enumerate}}
\def\beit{\begin{itemize}}
\def\enit{\end{itemize}}

\def\k_d{\kappa_d}

\def\co{\colon}

\def\k{\kappa}

\def\bepr{\begin{proposition}}
\def\enpr{\end{proposition}}
\def\vp{\varphi}

\def\Var{\ms{Var}}

\def\Tr{\ms{Tr}}

\usepackage{soul}
\def\Unif{\mathrm{Unif}}

\def\bN{\mathbb{N}}

\newcommand{\condPK}{\textsc{PK}}

\def\sech{\mathrm{sech}}

\begin{document}

\title{Central limit theorem for linear eigenvalue statistics of random spatial networks}

\author{Christian Hirsch}
\author{Kyeongsik Nam}
\author{Moritz Otto}
\author{Anne Marie Svane}

\begin{abstract}
The random connection model provides a flexible framework for spatial
networks in which connectivity is determined by both geometry and randomness.
Its vertices form a Poisson point process, and pairs are connected
independently, conditional on their positions, with distance-dependent
probabilities. This generalizes the hard distance threshold of random
geometric graphs and accommodates long-range connections. We also
consider a weighted extension in which each vertex carries
a positive random weight that influences its probability
of connecting to other vertices.
This combination of spatial structure and heterogeneous connectivity raises the question of whether spectral fluctuations exhibit Gaussian behavior on large scales.

We address this question by studying linear eigenvalue statistics $\operatorname{Tr}[f(A_n)]$ of the adjacency matrix $A_n$ of the random connection model in expanding windows $W_n$.
For polynomial test functions, we establish a finite
limit for the normalized variance and prove
a central limit theorem under suitable connection-profile
integrability and vertex mark moment assumptions.
When the limiting variance is positive, the normalized statistics admit  an
$O(|W_n|^{-1/2})$ normal-approximation bound in both
Wasserstein and Kolmogorov distance. Furthermore,
  under an exponential moment assumption on the vertex marks, we extend the central limit theorem to a general class of test functions $f$.
\end{abstract}

\maketitle

\section{Introduction}
\label{sec:intro}

Let $A$ be the adjacency matrix of a finite graph with eigenvalues
$\lambda_1,\ldots,\lambda_N$. For a test function $f$, the linear eigenvalue
statistic
\[
    \Tr[f(A)]=\sum_{i=1}^N f(\lambda_i)
\]
provides a basic summary of the graph spectrum. For polynomial $f$, these
statistics are combinations of closed-walk counts, while for more general
test functions they capture spectral information used in graph partitioning,
clustering, and related learning problems; see, for instance,
\cite{shi-malik,ng-jordan-weiss,von-luxburg}. Their fluctuations are well
understood for several classical random graph and random matrix models, where
independence or weak dependence between matrix entries allows the use of
methods from random matrix theory; see, for example,
\cite{furedi-komlos,erdos-knowles-yau-yin,shch1,cltsparse}. Much less is known
for sparse spatial graphs.  The shared random geometric feature introduces
nontrivial dependencies between edges.
These dependencies make the analysis of spectral
fluctuations considerably more challenging, particularly
in the sparse regime.

A fundamental example of a spatial network is the
random geometric graph, in which vertices are randomly
distributed and two vertices are
connected whenever their distance is below a fixed
threshold. We study its generalization, the \emph{random
connection model}, whose vertices form a Poisson point
process and whose edges are sampled independently,
conditional on the vertex positions, with probabilities
determined by a connection profile
$\varphi(\|x\|^d)$ for displacement $x$.
The hard-threshold profile recovers the random geometric
graph, whereas profiles with unbounded support permit
edges at arbitrarily large distances.  We also
consider a weighted extension in which every vertex carries a mark
and the marks of two endpoints modify the spatial scale of their connection. These
random geometric graphs and random connection models are among the standard
models of spatial networks; see \cite{penrose-rgg,barthelemy-spatial}.

The soft model presents a difficulty that is absent from the hard-threshold
setting. A random geometric graph has a deterministic interaction range, but
a soft random connection model generally does not: a vertex can be connected
to arbitrarily distant points. Moreover, inserting a deterministic vertex
requires generating all of its incident edges randomly while retaining the
same random edges among the original Poisson vertices. Thus, the local
representations commonly used for geometric Poisson functionals do not apply
directly. The question is whether the  spatial structure  still
 produces volume-order Gaussian fluctuations of the spectrum.

Our first result shows that it does. For polynomial test functions, we prove
that the variance of $\Tr[f(A_n)]$, divided by the observation volume, has a
finite limit. Under a stronger integrability assumption on the connection
profile and corresponding moment assumptions on the vertex marks, we then
prove a quantitative central limit theorem. When the limiting variance is
positive, the standardized statistic has Wasserstein and Kolmogorov distance
of order $|W_n|^{-1/2}$ from a standard Gaussian. The
integrability condition permits connection profiles with sufficiently fast
polynomial decay.

Our second objective is to move beyond polynomials. Under an exponential
moment assumption on the marks, we extend the central limit theorem to a
weighted Sobolev class of general test functions. The limiting variance is
obtained as the continuous extension of the variance form on polynomials.
Thus, the polynomial and general results describe two complementary aspects
of the same fluctuation theory: the first gives a quantitative normal
approximation under comparatively weak moment assumptions, while the second
gives a substantially larger class of spectral observables under stronger
control of the marks.

The most closely related result is the recent central limit theorem for
linear eigenvalue statistics of sparse random geometric graphs in
\cite{hirsch-nam-otto-rgg}. The present work replaces the deterministic
distance threshold by a soft connection profile and allows additional vertex
marks and interaction kernels. To the best of our knowledge, this gives the
first central limit theorems for linear eigenvalue statistics of sparse random
connection models with genuinely soft connection profiles. Existing spectral
results for spatial graphs have otherwise focused largely on limiting
spectral distributions or denser regimes; see, for example,
\cite{hamidouche-cottatellucci-avrachenkov,adhikari-adler-bobrowski-rosenthal}.
Spectral questions for geometric networks are also motivated by statistical
problems. Higher-order spectral clustering has been studied for soft
geometric block models in \cite{avrachenkov-bobu-dreveton}. This work
illustrates that spectral information in spatial networks is relevant not
only to global graph structure, but also to concrete inference problems.

The proof for polynomial test functions is based on two simple structural
ideas. First, we place all edge randomness on a common edge-marked probability
space. Insertion differences can then be formed without resampling edges
between the original Poisson vertices. Second, powers of the adjacency matrix
are expanded into closed walks. Products of insertion differences therefore
reduce to connected collections of finitely many walks, whose spatial
integrals can be controlled by passing to spanning trees and integrating their
vertices successively. This turns the dependence created by the geometry into
tree integrals involving only the connection profile and moments of the
marks. The resulting bounds allow us to apply the normal-approximation
estimates for edge-marked Poisson processes from \cite{nestmann}; related
stabilization methods for geometric Poisson and Gibbs functionals appear in
\cite{mal_stab,hirsch-otto-svane}.

The extension to general test functions uses a different mechanism. We
approximate the test function by polynomials and prove a uniform variance
estimate for exponentially twisted statistics. The main model-specific input
is an exponential trace bound. To obtain it, we orient the edges according to
the vertex marks and dominate the adjacency matrix, in quadratic-form order,
by a diagonal matrix involving the resulting incoming and outgoing degrees.
The exponential moment assumption on the marks then provides the required
degree estimates. Together with a Fourier--Duhamel argument, this yields the
variance bound needed to pass from polynomials to general test functions.

The remainder of the paper is organized as follows. Section~\ref{sec:results}
introduces the random connection model, its weighted extension, and the main
results. Section~\ref{subsec:variance_comparison} proves the variance
asymptotics and gives sufficient conditions for positivity of the limiting
variance. Section~\ref{sec:coupling_discretization} proves the quantitative
central limit theorem for polynomial test functions. Finally,
Section~\ref{sec:general} extends the central limit theorem to general test
functions.

\section{Model and Main Results}
\label{sec:results}

This section introduces the random connection model and its marked extension. It then defines the linear eigenvalue statistics and records the assumptions on the connection profile, interaction kernel, and marks. The section concludes with the variance asymptotics and central limit theorems that guide the later proofs.

%
%
\subsection{The random connection model}
\label{subsec:rcm}

The unmarked random connection model provides the basic geometric construction underlying the marked model used in the main results. Its vertices form a homogeneous Poisson point process in an expanding window, and distinct pairs are connected independently according to their distance. We first define this model and its adjacency matrix from \cite{netCom,pcPerc} before introducing the vertex marks and interaction kernel.

Fix \(d\ge1\). For \(n\ge1\), let \(\mathcal P_n\) be a homogeneous Poisson point process with intensity \(1\) on \(W_n:=[-n^{1/d}/2,n^{1/d}/2]^d\subseteq\mathbb R^d\). Then, \(|W_n|=n\), and the asymptotic regime is \(n\to\infty\) with the dimension \(d\) fixed. We denote by \(N_n:=|\mathcal P_n|\) the random number of points in \(\mathcal P_n\) and choose an arbitrary enumeration \(\mathcal P_n=\{X_1,\ldots,X_{N_n}\}\).

Let \(\varphi:[0,\infty)\to[0,1]\) be a non-increasing function, called the connection profile, satisfying 
\begin{align}\label{mdefinition}
    M_\vp:=\int_{\mathbb R^d}\varphi(\|x\|^d)\,\mathrm dx<\infty.
\end{align}  More generally, for \(\gamma\in(0,1]\), set
\begin{align}\label{mdef}
    M_{\varphi,\gamma}:=\int_{\mathbb R^d}\varphi(\|x\|^d)^\gamma\,\mathrm dx.
\end{align}

The random connection model (RCM) is the random graph with vertex set $\PP_n$. 
Conditionally on \(\mathcal P_n\), each pair of distinct points \(X_i,X_j\in\mathcal P_n\) is connected by an edge with probability \(\varphi(\|X_i-X_j\|^d)\), independently of all other pairs. The hard threshold case $\varphi(r)=\one\{r\le R\}$, for fixed $R>0$, produces the random geometric graph from \cite{hirsch-nam-otto-rgg} as a special case.
The adjacency matrix \(A_n\in\mathbb R^{N_n\times N_n}\) of the RCM is defined by
\begin{align}
A_n(i,j):=A_n(X_i,X_j):=\one\{\{X_i,X_j\}\text{ is an edge}\},\qquad 1\le i,j\le N_n,
\end{align}
with \(A_n(i,i)=0\). 

We work in an increasing-volume regime: the point intensity and
the connection profile remain fixed, while the observation window
expands. In particular, if a vertex is inserted at $x\in W_n$ with incident edges added independently with probabilities prescribed by $\varphi$, then the expected degree of $x$ is
\[
\int_{W_n}\varphi(\|x-y\|^d)\,dy
\le M_\varphi,
\]
uniformly in $n$ and $x$.
The expected number of edges is therefore at most
$M_\varphi|W_n|/2$, whereas the expected number of vertices is
$|W_n|$.
Thus, the RCM is sparse in the sense that the mean degree stays bounded as volume increases.

\subsection{The weighted RCM}

Defining the RCM as above, with the
hard-threshold random geometric graph as a special case, the degree of
a typical vertex follows a Poisson distribution.
This is a limitation when modeling networks with  degree
heterogeneity, in which highly connected hubs coexist with many
vertices of small degree. Power-law degree distributions are useful models
for such heterogeneity 
\cite{BroidoClauset2019}. 

Vertex weights provide a way to model differences in the tendency of individual vertices to form connections.
They may represent the attractiveness of a location
or the connectivity capacity of a communication node.
These interpretations motivate separating the spatial positions,
which encode geographical constraints, from the vertex weights.

We now precisely define the weighted RCM from \cite{gracarGrauerLuechtrathMoerters2019,komjathy2,gracarLuechtrathMoench2025}.
Attach to the points \(X_1,\ldots,X_{N_n}\) independent marks \(U_1,\ldots,U_{N_n}\) with common distribution \(\mu\) on \([1,\infty)\), independently of \(\mathcal P_n\). The resulting marked Poisson process is \(\check{\mathcal P}_n=\{(X_i,U_i):1\le i\le N_n\}\) on \(\check W_n:=W_n\times[1,\infty)\). We write \(\lambda(dx,du):=dx\,\mu(du)\) for its intensity measure. Throughout, \(U\) denotes a random variable with distribution \(\mu\).

Let \(\varphi\) be the connection profile introduced in Subsection~\ref{subsec:rcm}, and let \(\kappa:[1,\infty)^2\to[1,\infty)\) be a symmetric measurable function, called the interaction kernel. 
Given the marked configuration $\check{\mathcal P}_n$, the weighted RCM
has spatial vertex set $\mathcal P_n$.
For distinct marked points
$(X_i,U_i),(X_j,U_j)\in \check{\mathcal P}_n$, place an edge
between $X_i$ and $X_j$ with probability
\[
\varphi\!\left(
 \frac{\|X_i-X_j\|^d}{\kappa(U_i,U_j)}
\right),
\]
independently for different unordered pairs, conditional on
$\check P_n$. 

Intuitively, $\kappa(u,v)$ measures the tendency for two vertices with weights $u$ and $v$ to connect. For interaction kernels increasing in both arguments,
such as $\kappa(u,v)=uv$, vertices with larger weights have
a greater propensity to form connections and can interact
over larger spatial scales.
This mechanism underlies continuum scale-free percolation \cite{DeprezWuthrich2019}.

The unweighted RCM is recovered as a special case by taking \(U\equiv 1\)   and \(\kappa\equiv1\). 
The effect of the weights on the degree distribution can be seen
directly from the connection rule.
Consider the infinite-volume model and insert a vertex at the
origin with an independent mark $U_0$ having distribution $\mu $ and incident edges added independently.
Let $D$ be its degree.
For every $u$ such that $\mathbb E[\kappa(u,U)]<\infty$, the neighboring vertices form an 
independent thinning of the marked Poisson process and hence
\[
D\mid\{U_0=u\}
\sim\operatorname{Poisson}(\Lambda(u)),
\qquad
\Lambda(u)
:=
\int_{\mathbb R^d}\int_{[1,\infty)}
\varphi\!\left(\frac{\|x\|^d}{\kappa(u,v)}\right)
\,\mu (dv)\,dx
=
M_\varphi\,\mathbb E[\kappa(u,U)].
\]
Consequently, whenever $\Lambda(U_0)<\infty$ almost surely, the
unconditional degree distribution is mixed Poisson.
This provides a direct mechanism for introducing degree
heterogeneity while retaining a homogeneous spatial distribution
of vertices.

A particular  example is the kernel
$\kappa(u,v)=uv$.
If $0<M_\varphi<\infty$ and $\mathbb E U<\infty$, then
\[
\Lambda(u)=a u,
\qquad a:=M_\varphi\,\mathbb E U.
\]
Thus, the vertex weight is proportional to the conditional expected
degree.
For example, if the marks have a Pareto distribution,
$\mathbb P(U>t)=t^{-\beta}$ for $t\ge1$ with $\beta>1$, then
the degree distribution satisfies
\[
\mathbb P(D\ge k)
\sim
\mathbb P(aU\ge k)
=
a^\beta k^{-\beta},
\qquad k\to\infty.
\]
The degree therefore inherits the power-law tail of the marks. 
In particular, power-law degrees can arise from vertex
heterogeneity even when the connection profile itself has
bounded support.

~

Henceforth, $A_n$ denotes the adjacency matrix of this graph.
We identify a marked vertex $\check x=(x,u)$ with its spatial
coordinate when it is used as a vertex argument of a graph
quantity.
In particular,
\[
A_n(\check x,\check y):=A_n(x,y).
\]
For a graph with spatial vertex set $Q$ and its associated
vertex marks, write $\check Q$ for the corresponding marked
configuration. 

\subsection{Assumptions}
\label{subsec:assumptions}

The proofs require separate assumptions on the connection profile, the interaction kernel, and the mark distribution. We always assume that \eqref{mdefinition} holds for the connection profile. The quantitative CLTs further require that $M_{\varphi,1/4}<\infty$ with $M_{\varphi,1/4}$ defined in \eqref{mdef}. This stronger condition remains
compatible with long-range connections.
For example, for
\[
\varphi(t)=(1+t)^{-\alpha},\qquad \alpha>1,
\]
polar coordinate integration yields
\[
M_{\varphi,\gamma}<\infty
\quad\Longleftrightarrow\quad
\alpha\gamma>1.
\]
Thus, the condition $M_{\varphi,1/4}<\infty$ permits this profile
whenever $\alpha>4$.
The absence of $d$ from this threshold reflects our convention of
evaluating the profile at the $d$th power of the Euclidean distance.

For polynomial test functions, we assume finite moments of the vertex marks depending on the polynomial degree.
For general test functions, we impose the following stronger
exponential-moment condition, which is satisfied, for example,
by shifted exponential distributions.
We say that the mark distribution satisfies
Condition \eqref{cond:EM} if there exists \(\eta>0\) such that
\begin{equation}
\tag{EM}\label{cond:EM}
\mathbb E[e^{\eta U}]<\infty.
\end{equation}

The key assumption on our kernel $\kappa$ will be that it is dominated by the \emph{product kernel}. More precisely, we say that the interaction kernel satisfies Condition \eqref{cond:PK} if there
exists \(c_\kappa\ge1\) such that 
\begin{equation}
\tag{PK}\label{cond:PK}
1\le\kappa(u,v)\le c_\kappa uv
\qquad\text{for all }u,v\ge1.
\end{equation}
The assumptions separate two sources of variability.
Integrability of the connection profile controls the spatial extent
of interactions, whereas moments of the marks control the
heterogeneity of the interaction scales.
Condition \eqref{cond:PK} links these two controls: integrating one connection
probability over a spatial coordinate produces the factor
$M_\varphi\kappa(u,v)$, which is at most
$c_\kappa M_\varphi uv$.
Repeated spatial integration therefore introduces powers of the
vertex marks. This explains why the required mark moments depend
on the degree of the polynomial test function (see Theorem \ref{thm:polynomials} below).

These conditions are sufficient for the estimates used below;
we do not claim that they are necessary for Gaussian fluctuations.

\subsection{Linear eigenvalue statistics}

Spectrum  of
the adjacency  matrix $A_n$
is the random object studied throughout the paper. This subsection introduces the associated linear eigenvalue statistics and fixes the notation used in the main results. 

Let \(\lambda_1,\ldots,\lambda_{N_n}\) denote the real eigenvalues of \(A_n\). For a test function \(f:\mathbb R\to\mathbb R\), define
\begin{align}
L_n(f):=\sum_{i=1}^{N_n}f(\lambda_i)=\Tr[f(A_n)].
\end{align}
Our goal is to establish a finite limit for the variance of $L_n(f)$
divided by $|W_n|$ and a central limit theorem for the
centered statistic normalized by $\sqrt{|W_n|}$.  When the limiting variance is positive, this also yields a standard normal approximation under normalization by the finite-volume standard deviation. 

\subsection{Main theorems}

This subsection states the three principal results of the paper. The first establishes variance asymptotics for polynomial test functions, the second proves a quantitative central limit theorem in that class, and the third extends qualitative convergence to a general test-function space. Their order reflects the proof structure: the variance limit enters the polynomial central limit theorem, which then provides the approximations used for general test functions.

Our first result establishes the existence of the limiting variance for polynomial test functions.

%
%
\begin{theorem}[Variance asymptotics]
\label{thm:variance_limit}
Let \(f \) be a polynomial of degree \(m\ge1\). Assume that \eqref{cond:PK} holds and that \(\mathbb E[U^{3m-1}]<\infty\). Then, the following limit exists and is finite:
\[
\lim_{n\to\infty}|W_n|^{-1}\Var(\Tr[f(A_n)])  .
\]
\end{theorem}
For every polynomial satisfying the hypotheses of Theorem~\ref{thm:variance_limit}, we denote this limit by \(\sigma_f^2\). Proposition \ref{prop:variance_existence_positivity} below provides conditions for positivity of $\sigma_f^2$.

Our second result establishes a quantitative central limit theorem for polynomial test functions. For integrable real-valued random variables \(X\) and \(Y\), let
\[
d_W(X,Y):=\sup_{\operatorname{Lip}(h)\le1}|\mathbb E[h(X)]-\mathbb E[h(Y)]|
\]
denote the Wasserstein distance between their distributions and let
\[
d_K(X,Y):=\sup_{t\in \R}|\P(X\le t)-\P(Y\le t)|
\]
denote their Kolmogorov distance. In addition, $\mathcal N(\mu, \sigma^2)$ denotes the normal distribution with mean $\mu$ and variance $\sigma^2.$
%
%
\begin{theorem}[CLT for polynomial test functions]
\label{thm:polynomials}
Let \(f\) be a polynomial of degree \(m\ge1\). Assume that \eqref{cond:PK} holds, that \(M_{\varphi,1/4}<\infty\), and that \(\mathbb E[U^{6m}]<\infty\). Then,
\[
\frac{\Tr[f(A_n)]-\mathbb E[\Tr[f(A_n)]]}{\sqrt{|W_n|}}
\xrightarrow{d}
\mathcal N(0,\sigma_f^2).
\]
If \(\sigma_f^2>0\), then there exists  a constant  \(C<\infty\) such that 
{\[
\max\left\{ d_W\left(
\frac{\Tr[f(A_n)]-\mathbb E[\Tr[f(A_n)]]}{\sqrt{\Var(\Tr[f(A_n)])}},
\mathcal N(0,1)
\right),d_K\left(
\frac{\Tr[f(A_n)]-\mathbb E[\Tr[f(A_n)]]}{\sqrt{\Var(\Tr[f(A_n)])}},
\mathcal N(0,1)
\right) \right\}
\le \frac{C}{\sqrt{|W_n|}},
\]}
where \(C\) may depend on \(f\) and the fixed model parameters, but not on \(n\).
\end{theorem}

Our third result extends the central limit theorem to a broader class of test functions. 
For \(c>0\), let \(\mathcal H_c^2\) be the space of twice weakly differentiable functions \(f:\mathbb R\to\mathbb R\) such that
\[
\|f\|_{\mathcal H_c^2}^2
:=
\int_{\mathbb R}(|f(x)|^2+|f'(x)|^2+|f''(x)|^2)\sech^2(cx)\,\mathrm dx
<\infty.
\] 
We further set \(\mathcal H_{c,0}^2:=\{f\in\mathcal H_c^2:f(0)=0\}\). {Every element of $\mathcal  H_c^2$ belongs to
$H^2_{\mathrm{loc}}(\mathbb R)$, since the weight
$\operatorname{sech}^2(cx)$ is bounded below on compact intervals.
We always use its unique locally $C^1$ representative.
In particular, $f(0)$ and the spectral values
$f(\lambda_i(A_n))$ refer to this representative.}

The exponential weight serves two purposes.
First, unlike the usual unweighted Sobolev space,
$\mathcal H_c^2$ contains every polynomial, so it provides a common space
for the polynomial CLT and its extension.
Second, it allows test functions that grow at infinity, while
still controlling the two derivatives.
The condition $f(0)=0$ is convenient when comparing adjacency
matrices before and after a vertex insertion.
Adding an isolated vertex appends a zero eigenvalue, and this
condition ensures that such insertion does not change the statistic. 
The following theorem extends the polynomial CLT to this
test-function space.

%
%
\begin{theorem}[CLT for general test functions]
\label{thm:general_functions}
Assume that \eqref{cond:EM} and \eqref{cond:PK} hold and that \(M_{\varphi,1/4}<\infty\). Let \(c>0\) and \(f\in\mathcal H_{c,0}^2\).   Then, there exists \(\sigma_f^2\ge0\) such that
\[
\frac{\Tr[f(A_n)]-\mathbb E[\Tr[f(A_n)]]}{\sqrt{|W_n|}}
\xrightarrow{d}
\mathcal N(0,\sigma_f^2)
\qquad\text{as }n\to\infty.
\]
\end{theorem}

The proof identifies \(\sigma_f^2\) with the value at \(f\) of the continuous extension of the limiting variance form on polynomials.

\subsection{Acknowledgments}
C.H.~is supported by a research grant (VIL69126) from VILLUM FONDEN. K.N. is supported by Samsung Science and Technology Foundation under
Project Number SSTF-BA2202-02.  M.O. is supported by the NWO Gravitation project NETWORKS under grant agreement No. 024.002.003. 
This
article is based upon work from COST Action 24122 mSPACE, supported by COST (European Cooperation in Science and Technology).

%
%
\section{Proof of the variance asymptotics}
\label{subsec:variance_comparison}

This section proves the existence of a limiting volume-normalized variance and records a sufficient criterion for its positivity. The proof relies on variance asymptotics results for weakly stabilizing functionals from  \cite{trinh2, pabst2026betti}. Verifying the conditions of these results involves studying the difference when removing a fixed vertex from the RCM. To make this rigorous, we recall a construction from \cite{nestmann} of the RCM and its difference operators from an edge-marked Poisson process in Section \ref{ss:diff}. This construction will be used again for the proof of the CLT in Theorem \ref{thm:polynomials}.  We recall the main result from \cite{trinh2,pabst2026betti} in Section \ref{ss:weak} before proving the main theorem in Section \ref{ss:variance_proof}.

\subsection{Notation}\label{ss:notation}
Let $G$ denote a finite graph with (deterministic) vertex set $P$ and adjacency matrix $A$. We denote the entry of $A$ corresponding to $x_1,x_2\in P$ by $A(x_1,x_2;G)$. Consider the linear eigenvalue statistic $L(f;G) = \Tr[f(A)]$ with a polynomial test function $f$.
For integers $1\le p\le m$, define
\begin{align} 
  \Sigma_{m,p} := \big\{\pi\co \{1, \dots, m\} \to \{1, \dots, p\}\text{ surjective}\big\}.
\end{align}
This corresponds to the set of ordered partitions of $\{1, \dots, m\}$ into $p$ nonempty sets.
We also define
\[
\Sigma^\circ_{m,p}
:=\bigl\{\pi\in\Sigma_{m,p}:
          \pi(j)\ne\pi(j+1)\text{ for }j=1,\ldots,m\bigr\},
\]
where we use the convention $\pi(m+1)  = \pi(1)$.
Let $x_1,\dots,x_p$ be $p$ mutually distinct vertices in $G$ and set 
\begin{align*}
	g_{m, p}(x_1,\dots, x_p;G):={\f1{p!}} \sum_{\pi \in \Sigma_{m,p}^\circ}\prod_{j =1}^{m} A(x_{\pi(j)}, x_{\pi(j + 1)} ;G).
\end{align*}
 Thus, $ p! g_{m, p}(x_1,\dots, x_p)$ counts based
closed walks of length $m$ whose visited vertex set
is exactly $\{x_1,\ldots,x_p\}$.

For $m \ge 1$, one can write $L(x^m;G)=\Tr[A^m]$ as a sum of $m - 1$ U-statistics
\begin{align*}
	L(x^m;G) &= \sum_{p =2}^ m \sum_{(X_1,\dots,X_p) \in P^{(p)}}\f1{p!}\sum_{\pi \in \Sigma_{m,p}^\circ}\prod_{j =1}^ mA(X_{\pi(j)}, X_{\pi(j + 1)};G) 
  = \sum_{p =2}^ m \sum_{(X_1,\dots,X_p) \in P^{(p)}}g_{m, p}(X_1,\dots, X_p; G), 
\end{align*}
where $P^{(p)}$ denotes the set of $p$-tuples of pairwise distinct points in $P$. Thus, $L(x^m;G)$ counts the total number of closed walks of length $m$ in the graph.  To see why this holds,  
start from the matrix-product identity
\[
L(x^m;G) =\operatorname{Tr}[A^m]
=
\sum_{(v_1,\ldots,v_m)\in P^m}
\prod_{j=1}^m A(v_j,v_{j+1};G),
\qquad v_{m+1}:=v_1.
\]
Group the vertex sequences $(v_1,\ldots,v_m)$ according to the number $p$
of distinct vertices they visit, where
 $p\ge2$ because diagonal entries of
$A$ vanish.
For a fixed sequence visiting exactly $p$ vertices, choose
an ordered enumeration $(X_1,\ldots,X_p)$ of those vertices.
Then there is a unique surjection
$\pi:\{1,\ldots,m\}\to\{1,\ldots,p\}$ such that
$v_j=X_{\pi(j)}$ for every $j$. As   diagonal entries of the
adjacency matrix  are zero, only terms with $\pi\in \Sigma^\circ_{m,p}$ are non-vanishing.
Conversely, the ordered tuple and the surjection recover
the original vertex sequence.
There are exactly $p!$ possible ordered enumerations of its
visited vertex set. 
The factor $1/p!$ removes  this multiplicity.

By linearity of $L(\cdot;G)$, a polynomial $f(x)=a_0+\cdots+a_mx^m$ satisfies 
\begin{align}\nonumber
 L(f;G):=\Tr[f(A)]& = a_0 |P| + \sum_{q = 2}^m a_q \sum_{p = 2}^q\sum_{(X_1,\dots,X_p) \in P^{(p)}} g_{q, p}(X_1,\dots, X_p;G) \nonumber \\
 &= a_0 |P| + \sum_{p = 2}^m\sum_{(X_1,\dots,X_p) \in P^{(p)}}  \sum_{q = p}^ma_q  g_{q, p} (X_1,\dots, X_p;G).\label{eq:sumscore}
\end{align}

\subsection{Construction of the RCM from edge marked sequences} \label{ss:diff}
We first introduce a construction of the RCM from \cite{nestmann} that allows a convenient way of adding extra fixed points to the graph.
Define a measurable total ordering $\preceq$ on $\R^d$ as follows: order two points by their $\ell_\infty$-distance to the origin. Two points with the same distance to the origin are ordered lexicographically. Let $\check \PP$ be a Poisson process on $\check \R^d = \R^d \times [1,\infty)  $ with intensity measure $\la$ and let $\PP$ denote the underlying Poisson process on $\R^d$. The ordering on $\R^d$ almost surely allows us to label the points  in $\PP$ as $\PP=\{X_i, \,i\in \N\}$ with $X_i\prec X_j$ for $i<j$.  Let $[0,1]^{\Z\times \Z}$ denote the space of doubly indexed sequences equipped with the product of uniform measures such that a random element $\mathbf R=(R_{i,j})_{i,j\in\mathbb Z}$ in $[0,1]^{\Z\times \Z}$ corresponds to a doubly indexed sequence of i.i.d.\ uniform random variables $R_{i,j}$. Let $\mathbf R$ be such a random element which is independent of the entire
marked Poisson process $\check \PP$. Define the edge-marked Poisson process as the pair
$$\xi= (\check \PP, \mathbf R). $$
We refer to $\mathbf R$ as the edge-mark sequence; each mark $R_{i,j}$ with $i<j$ is associated with the pair of Poisson points $X_i$ and $X_j$.
We can build a weighted random connection model $G(\xi)$ as the graph with vertex set $\PP$. If $(X_i,U_i),(X_j,U_j)\in \check \PP$ with $i<j$, we put an edge between $X_i$ and $X_j$ if 
\begin{align} \label{conrule}
    R_{i,j}\le \vp\Big(\frac{\|X_i-X_j\|^d}{\kappa(U_i,U_j)}\Big).
\end{align}
Fix an ordered tuple
\[
\check  \xx_r:=(\check x_1,\ldots,\check x_r),
\qquad
\check x_i=(x_i,u_i)\in\mathbb R^d\times[1,\infty),
\]
with pairwise distinct spatial coordinates.
The tuple order determines negative labels
$-1,\ldots,-r$; that is, we write  $\check x_i =(x_{i},u_i)=(X_{-i},U_{-i})$. 
For $I\subseteq \{1,\ldots,r\}$, we construct a random graph $G_I(\xi, \check \xx_r)$ as follows:  The vertex set of $G_I(\xi,\check \xx_r)$ is $\PP\cup \{X_{-i}:i\in I\}$.  For
$i,j\in\mathbb N\cup\{-a:a\in I\}$ with $i<j$,
we add an edge between $X_i$ and $X_j$ when   \eqref{conrule} holds.
This construction ensures that if $J\subseteq I$, then the graph $G_I(\xi, \check \xx_r)\setminus \{x_i,i\in J\}$,  where we remove the deterministic vertices indexed by $J$ and all incident edges from $G_I(\xi, \check \xx_r)$, equals $G_{I\setminus J}(\xi, \check \xx_r)$. In particular, $G_I(\xi, \check \xx_r)\setminus \{x_i,i\in I\}=G_\emptyset(\xi)=G(\xi)$; that is,  we recover the  RCM $G(\xi)$ if we remove all the fixed vertices in $\check \xx_r$.

This coupling of RCMs with one or more fixed points added allows a notion of difference operators adapted to the RCM setting that is useful to establish asymptotic results \cite{nestmann}. Let $H$ be a real-valued measurable functional on  graphs
with vertices in $\R^d$, and set $F:=H(G(\xi)).$  
Define the difference operators
\[
\Delta^r_{\check x_1,\ldots,\check x_r; \check\xx_r}F
=
\sum_{I\subseteq\{1,\ldots,r\}}
(-1)^{r-|I|}
H\bigl(G_I(\xi,\check \xx_r)\bigr).
\] 
We also consider differences with respect to a subset
of the inserted vertices while retaining the full
ambient tuple $\check \xx_r$ to specify
their reserved labels.
For $i,j\in \{1,\ldots, r\}$ with $i\neq j$, we thus write
\begin{align*}
&\Delta_{\check x_i;\check \xx_r}F=\Delta^1_{\check x_i;\check \xx_r}F
= H\bigl(G_{\{i\}}(\xi,\check \xx_r)\bigr)-H\bigl(G_{\emptyset}(\xi,\check \xx_r)\bigr),\\
&\Delta^2_{\check x_i,\check x_j;\check \xx_r}F
= H\bigl(G_{\{i,j\}}(\xi,\check \xx_r)\bigr)-H\bigl(G_{\{i\}}(\xi,\check \xx_r)\bigr)- H\bigl(G_{\{j\}}(\xi,\check \xx_r)\bigr)+H\bigl(G_{\emptyset}(\xi,\check \xx_r)\bigr).
\end{align*}

For  $W\subseteq \R^d$,  we let $\check W = W\times [1,\infty)$ and define $\PP_W=\PP\cap W$, $\check \PP_W=\check \PP\cap \check W$, and $\xi_W=(\check \PP_W,\mathbf R)$. Let $G_{W,I}(\xi,\check \xx_r)$ denote the restriction of $G_I(\xi,\check \xx_r)$ to the vertex set $(\mathcal P\cup\{x_i:i\in I\})\cap W$.  When $W=W_n$, we simply write $\PP_n$, $\xi_n$, $G_n$, $G_{n,I}$ etc. Note that $G_{n,I}(\xi,\check \xx_r)=G_I(\xi_n,\check \xx_r)$ when $\check{\xx}_r\in\check W_n^r$. Often, we will consider a functional $F$ that depends on $G(\xi)$ only via $G_W(\xi)$. In this case, we may think of the restriction as part of the functional $H$.

\medskip
\noindent\textbf{Multivariate Mecke formula.}
Fix an ordered tuple
$\check \xx_r=(\check x_1,\ldots,\check x_r)$
of marked vertices in $\mathbb R^d\times[1,\infty)$
with distinct spatial coordinates. Write
\[
G_{n,\check{\xx}_r}
:=
G_{n,\{1,\ldots,r\}}(\xi,\check{\xx}_r),
\qquad
G_{n,\check{\xx}_0}:=G(\xi_n).
\] 
  For any nonnegative measurable graph functional 
$h(\check{\yy}_k;G)$, invariant under
auxiliary relabeling, 
\begin{align} \label{eq:rcm-mecke}
&\mathbb E \Big[ 
\sum_{\check{\yy}_k\in\check\PP_n^{(k)}}
h\bigl(
 \check{\yy}_k;
 G_{n,\check{\xx}_r}
\bigr)  \Big]=
\int_{\check W_n^k}
\mathbb E\left[
 h\bigl(
  \check{\yy}_k;
 G_{n,\check{\xx}_r,\check{\yy}_k}
 \bigr)
\right]
\prod_{j=1}^k\lambda(d\check y_j).
\end{align}
Here $G_{n,\check{\xx}_r,\check{\yy}_k}$ denotes the graph obtained
by inserting the $r+k$ marked vertices in the concatenated tuple
$(\check{\xx}_r,\check{\yy}_k)$ into $G(\xi_n)$ (the order of the concatenated tuple specifies the 
labels, but does not affect
the conditional distribution of the resulting graph).
The sum is over ordered $k$-tuples of distinct Poisson vertices,
while the inserted vertices $\check{\xx}_r$ remain fixed.
The case $r=0$ corresponds to the original RCM.

To justify the identity, first condition on the
marked vertex configuration and average the edge
randomness.
The resulting function depends only on the marked
vertices, so the ordinary multivariate Mecke
formula \cite[Theorem 4.4]{last2018lectures} applies.
Conditional on the augmented marked configuration,
the graph  has  
independent edge distributions, giving the identity.

\subsection{A result from the literature}\label{ss:weak}
Let $W\subseteq\mathbb R^d$ be a bounded Borel set.
If the origin $o$ belongs to $W$ and $u\in[1,\infty)$, then $G_{W,\{1\}}(\xi,((o,u)))$ is the
weighted RCM restricted to $W$ with the additional
vertex $(o,u)$ and its incident edges.
If $o\notin W$, the inserted vertex is removed by the
spatial restriction, so the augmented and unaugmented
restricted graphs coincide.
For a polynomial $f$, define the add-one cost 
\begin{equation}\label{eq:add-one}
  \Delta_{W,f}(u):=\Delta^1_{(o,u);((o,u))} L(f;G_W(\xi)) = L(f; G_{W,\{1\}}(\xi,((o,u)))) - L(f; G_{W,\emptyset}(\xi,((o,u)))).
\end{equation}
Let $\mathcal{W}$ denote the collection of cubes of the form $\prod_{i=1}^d [x_i,x_i+l)$, where $x_1,\ldots,x_d\in\R$ and $l>0$.  

The main ingredient in the proof of Theorem \ref{thm:variance_limit} is the following result. It is 
essentially the variance convergence part of the weak stabilization CLT \cite[Thm. 1.1]{trinh2} in the unweighted setting, while the weighted case is covered by \cite[Lem. 2.2]{pabst2026betti}. We state it here in a version adapted to the functional $L(f;G)$. 

\begin{theorem}[Variance asymptotics from \cite{trinh2,pabst2026betti}]\label{thm:trinh}
Let $f$ be a polynomial test function.
Suppose that 
\begin{itemize}
\item[(i)] The limit $\Delta_f(u):=\lim_{W^n\uparrow \R^d}\Delta_{W^n,f}(u)$ exists almost surely for any sequence of cubes $W^n\in \mathcal{W}$ with $W^n\uparrow \R^d$ and any mark $u$.
\item[(ii)] For some $p>2$, $\sup_{W\in \mathcal{W}} \E[|\Delta_{W,f}(U)|^p]<\infty$, where $U$ has law $\mu$ and is independent of
the original marked Poisson process and the
edge randomness.
\end{itemize}
If (i) and (ii) hold, then the limit
$$\sigma_f^2:=\lim_{n\to\infty}|W_n|^{-1}\Var(L(f;G_n(\xi))),$$ 
exists, and 
$\sigma_f^2>0$ whenever $\mathbb P(\Delta_f(U)\ne0)>0$.
\end{theorem}

\begin{proof}
The polynomial linear eigenvalue statistic $L(f;G)$ is invariant under translations of the vertex set of $G$. Hence, the variance asymptotics part of \cite[Thm. 1.1]{trinh2} gives the result in the unweighted case. The weighted case follows from \cite{pabst2026betti}, where the RCM is the special case $\alpha=1$ of the random simplicial complexes introduced there. To see how the result follows, note that $G_{\{1\}}(\xi,(o,u))$ has the same distribution as the random graph called $\Delta^{(o,u,1,V)}$ in \cite{pabst2026betti}. Here $ V$ denotes the independent auxiliary mark used to sample
the edges. It follows that the family of random variables  $\{\Lambda_{(o,u,1,V)}L(f,\Delta_W)\}_{W\in \mathcal{W}}$ in \cite{pabst2026betti} have the same joint distribution as our $\{\Delta_{W,f}(u)\}_{W\in \mathcal{W}}$ and thus Theorem \ref{thm:trinh} follows directly from \cite[Lem. 2.2]{pabst2026betti}.
\end{proof}

\subsection{Proof of Theorem \ref{thm:variance_limit}}\label{ss:variance_proof}

We show Theorem \ref{thm:variance_limit} by showing  the following proposition.

\begin{proposition}[Variance asymptotics and positivity]
\label{prop:variance_existence_positivity}
Let \(f(x)=\sum_{q=0}^m a_qx^q\) be a polynomial of degree \(m\ge1\). Assume that \eqref{cond:PK} holds and that \(\mathbb E[U^{3m-1}]<\infty\). Then, the limiting add-one cost $\Delta_f(u)$ exists almost surely for every $u\ge1$, and the limit $\lim_{n\to\infty}|W_n|^{-1}\Var(\Tr[f(A_n)])$ exists.

If $\P(\Delta_f(U)\neq0)>0$, where $U$ has law $\mu$, then $\sigma_f^2>0$. If $M_\vp>0$ and either $a_0\ne0$ or all coefficients of $f$ are non-negative or all are non-positive, and $f$ is not of the form $f(x)=a_1x$, then $\sigma_f^2>0$.
\end{proposition}

\begin{proof}[Proof of Theorem~\ref{thm:variance_limit}]
The first conclusion of Proposition~\ref{prop:variance_existence_positivity} gives the existence and finiteness of the limit asserted in Theorem \ref{thm:variance_limit} under exactly the same hypotheses.
\end{proof}

The proof of Proposition \ref{prop:variance_existence_positivity} represents polynomial spectral statistics by closed-walk counts as in Section \ref{ss:notation}. Their moments are controlled using the following tree estimate, which relies on an idea of \cite[Theorem 2.1]{heerten-hirsch-otto}.

\begin{lemma}[Tree integral estimate]
\label{lem:weighted_tree_estimate}
Let \(\gamma\in(0,1]\), assume that \(\kappa\) satisfies \eqref{cond:PK}, and suppose that \(M_{\varphi,\gamma}<\infty\). Let \(\mathcal T\) be a   tree with vertex set \(\{0,1,\ldots,k\}\), and let \(d_i\) denote the degree of vertex \(i\) in \(\mathcal T\). If \(u_0,\ldots,u_k\in[1,\infty)\) and \(x_0\in\mathbb R^d\), then
\[
    \int_{(\mathbb R^d)^k}
    \prod_{\{i,j\}\in E(\mathcal T)}
    \varphi\bigg(\frac{\|x_i-x_j\|^d}{\kappa(u_i,u_j)}\bigg)^\gamma
    \,\mathrm dx_1\cdots\mathrm dx_k
    \le
    M_{\varphi,\gamma}^k c_\kappa^k\prod_{i=0}^k u_i^{d_i}.
\]
\end{lemma}

\begin{proof}[Proof of Proposition \ref{prop:variance_existence_positivity}]
 For the first statement, we need to show (i) and (ii) in Theorem \ref{thm:trinh}. If the origin $o$ is not contained in $W$, then $\Delta_{W,f}(u)$ is zero.
It is therefore enough to prove the uniform moment
bounds for cubes containing the origin.
For such cubes, linearity gives   
    \begin{equation}\label{eq:lin-Delta}
    \Delta_{W,f}(u) = a_0 + a_2\Delta_{W,2}(u) + \ldots +a_m\Delta_{W,m}(u)
    \end{equation}
    where $\Delta_{W,q}(u):=\Delta_{W,x^q}(u)$. Here we used that $\Delta_{W,0}(u)=1$, and $\Delta_{W,1}(u)=0$ because the graph has no loops. It follows that it is enough to show  (i) and (ii)  in Theorem \ref{thm:trinh} for $f(x)=x^q$ where $ 2 \le q \le m$.

    First consider (i). The loop counting interpretation of the functional $L(x^q,\cdot)$ from Section \ref{ss:notation} implies that $\Delta_{W,q}(u)$ counts the number of closed walks of length $q$ in $G_{W,\{1\}}(\xi ,(o,u))$ that pass through $o$ at least once. By multivariate Mecke formula, the expected number of rooted closed walks of length $q$ starting and ending at $o$ is
    \begin{align} \nonumber
&\E\bigg[\sum_{p=2}^{q}\sum_{((X_2,U_2),\dots,(X_p,U_p))\in(\check\PP\cap\check W)^{(p-1)}}\f1{(p-1)!}\sum_{\substack{\pi\in\Sigma_{q,p}^\circ\\\pi(1)=1}}\prod_{j=1}^q A(X_{\pi(j)},X_{\pi(j+1)};G_{W,\{1\}}(\xi, (o,u) )) \bigg]\\
    &=\sum_{p=2}^{q}\f1{(p-1)!}\sum_{\substack{\pi\in\Sigma_{q,p}^\circ\\\pi(1)=1}}\int_{W^{p-1}}\E \bigg[ \prod_{\{i,j\}\in E(\pi)}\vp\bigg(\frac{\|x_i-x_j\|^d}{\kappa(U_i,U_j)}\bigg) \bigg]\,\d x_2\cdots\d x_p,\label{eq:allpaths}
    \end{align}
where $X_1=x_1=o$, $U_1=u$, and $E(\pi)$ is the set of distinct undirected edges traversed by the closed walk. The graph with vertex set $\{1,\ldots,p\}$ and edge set $E(\pi)$ is connected. Choose a spanning tree $\mathcal T$ of this graph with edge set $E(\mathcal{T})\subseteq E(\pi)$. Since every connection probability lies in $[0,1]$,
\begin{align}
&\int_{W^{p-1}}\prod_{\{i,j\}\in E(\pi)}\vp\bigg(\frac{\|x_i-x_j\|^d}{\kappa(U_i,U_j)}\bigg)\,\d x_2\cdots\d x_p
\le\int_{W^{p-1}}\prod_{\{i,j\}\in E(\mathcal T)}\vp\bigg(\frac{\|x_i-x_j\|^d}{\kappa(U_i,U_j)}\bigg)\,\d x_2\cdots\d x_p.\label{eq:onepath}
\end{align} 
 By Lemma \ref{lem:weighted_tree_estimate}, the expectation of \eqref{eq:onepath} is bounded by
\(
 u^{p-1}{M_\vp^{p-1}c_\kappa^{p-1}} \prod_{i=2}^{p}\E[U_i^{d_i}],
\)
where $ d_i \le p-1$. 
Let $C_{q,W}(u)$ denote the number of length-$q$
closed walks based at $(o,u)$ in the common augmented
graph restricted to $W$, where $o\in W$.
Its expectation is the expression in \eqref{eq:allpaths}.
The preceding  bound is independent of $W$,
and the number of terms in \eqref{eq:allpaths} depends only on $q$.
Since all required mark moments have order at most
$q-1\le3m-1$, we have $\sup_n\mathbb E C_{q,W_n}(u)
\le C_q u^{q-1}<\infty $     for each  $u\ge1$.

Let $C_{q,\infty}(u)$ count the corresponding rooted
walks in the infinite augmented graph.
All finite-window graphs are restrictions of this
same graph, and every finite walk is contained
in a sufficiently large window.
Consequently, $C_{q,W_n}(u)\uparrow C_{q,\infty}(u).$ 
By monotone convergence,
\[
\mathbb E C_{q,\infty}(u)
=
\lim_{n\to\infty}\mathbb E C_{q,W_n}(u)
\le C_q u^{q-1}<\infty.
\]
Thus $C_{q,\infty}(u)<\infty$ almost surely.
The union of the vertices visited by these finitely
many rooted walks is finite.
Every based closed walk of length $q$ visiting
the inserted vertex is a cyclic shift of a walk
based there, so it visits vertices in the same
finite union.
Therefore, this union contains every vertex that
can contribute to the add-one cost and is contained
in some random window $W_N$. 
It follows that
$
  \Delta_{W,q}(u)=\Delta_{W_N,q}(u)
$
whenever $W_N\subseteq W$. Hence,
\[
    \Delta_q(u):=\lim_{W^n\uparrow\mathbb R^d}\Delta_{W^n,q}(u) = \Delta_{W_N,q}(u)
\]
exists almost surely, which shows (i).

~

    For the moment condition (ii), we show that  
    $ \sup_{W\in \mathcal{W}} \E[|\Delta_{W,q}(U)|^3]<\infty$.   We have
 \begin{align} 
 &\E[|\Delta_{W,q}(U)|^3] \nonumber \\
 &\le \E \bigg[ q^3\bigg(\sum_{p =2}^{q} \sum_{((X_2,U_2),\dots,(X_p,U_p)) \in (\check \PP\cap \check W)^{(p-1)}}\frac{1}{(p-1)!}\sum_{\substack{\pi \in \Sigma_{q,p}^\circ\\ \pi(1)=1}}\prod_{j =1}^q A(X_{\pi(j)}, X_{\pi(j + 1)};G_{W,\{1\}}(\xi,(o,U)) ) \bigg)^3\bigg],\label{eq:3paths}
 \end{align}
  where the factor $q^3$ accounts for the at most $q$ possible starting points of the same closed walk.
   Expanding \eqref{eq:3paths}, we obtain a sum in which each term corresponds to three rooted closed walks of length $q$, each starting and ending at $o$. Concatenating them at $o$ produces a closed walk of length $3q$ that starts and ends at $o$. Thus, \eqref{eq:3paths} is bounded by $q^3$ times the expected total number of rooted closed walks of length $3q$ starting and ending at $o$. For fixed $u$, this was bounded by \eqref{eq:allpaths} with $q$ replaced by $3q$, and each integral was again bounded. After concatenation, the resulting walk has length $3q$, so the weighted tree estimate involves only moments of order at most $3q-1\le3m-1$. Recalling the condition $\mathbb E[U^{3m-1}]<\infty$, taking expectation with respect to $U$ in \eqref{eq:allpaths} as above therefore yields a finite bound independent of $W$.  This shows (ii).

~

It remains to verify the conditions for variance positivity. Clearly, $\Delta_0=1$ and $\Delta_1=0$. The probability that there are no edges incident to $(o,u)$, for fixed $u$, is
\begin{align*}&\E\bigg[\prod_{(X_i,U_i)\in \check \PP} (1- A(o,X_i))\bigg] = \E\bigg[\prod_{(X_i,U_i)\in \check \PP} \Big(1- \vp\Big(\frac{\|X_i\|^d}{\kappa(u,U_i)}\Big)\Big)\bigg]\\
&= \exp\Big( - \int_{\R^d} \E \Big[\vp\Big(\frac{\|x\|^d}{\kappa(u,U)}\Big)\Big]\d x\Big)\ge \exp(-c_\kappa M_\vp u\E[U])
\end{align*}
by Lemma \ref{lem:weighted_tree_estimate}.
Since  $M_\vp <\infty$ by the assumption \eqref{mdefinition},  there is a positive probability that no edges are incident to $(o,u)$, in which case $\Delta_q(u)=0$ for all $q>0$, and hence $\P(\Delta_f(U) =a_0)>0$. If $a_0\neq 0$, this implies $\sigma_f^2>0$.

For the second case, note that $\Delta_q(u)\ge0$ for every $q\ge2$, because adding a vertex at $(o,u)$ together with its incident edges does not destroy existing closed walks but may add extra closed walks that pass through $o$. 
A walk of length 3 starting and ending at $o$ must pass through exactly two other points; hence the expected number of such paths is
\begin{align*}
 \E &\bigg[ \sum_{((X_1,U_1),(X_2,U_2))\in \check \PP^{(2)}} A(o,X_1)A(X_1,X_2)A(X_2,o) \bigg] \ge \int_{\R^d}\int_{\R^d} \vp(\|x_1\|^d)\vp(\|x_1-x_2\|^d)\vp(\|x_2\|^d)\d x_1\d x_2 ,\label{eq:3cycle}
\end{align*}
 where we used that $\vp$ is nonincreasing and $\kappa(u,v)\ge 1$.
Since $\vp$ is nonincreasing and $M_\vp>0$, $\vp$ must be positive on a ball of (possibly infinite) radius $\rho>0$. It follows that the integrand above must be positive on $B_{\rho/2}(o)\times B_{\rho/2}(o)$ and hence the integral is positive. Thus, $\P(\Delta_3(U)>0)>0$. Going back and forth along an edge incident to $o$ creates a closed walk of length 2. For any $q>3$, a closed walk of length $q$ can be constructed by concatenating closed walks of lengths $2$ and $3$. Therefore, $\Delta_3(u)>0$ implies $\Delta_q(u)> 0$ for all $q> 1$.
It follows from \eqref{eq:lin-Delta} that if all coefficients are non-negative (resp. non-positive) and $f$ is not of the form $a_1x$, then $\Delta_f(u)$ is also positive (resp. negative) on the event $\Delta_3(u)>0$.

\end{proof}

We finally establish Lemma~\ref{lem:weighted_tree_estimate}.

\begin{proof}[Proof of Lemma~\ref{lem:weighted_tree_estimate}]
We integrate the spatial variables successively, beginning with a leaf of $\mathcal T$ other than the fixed vertex $0$. If $j$ is a leaf with parent $i$, then the change of variables $z=\kappa(u_i,u_j)^{-1/d}(x_j-x_i)$ gives
\[
    \int_{\mathbb R^d}
    \varphi\bigg(\frac{\|x_j-x_i\|^d}{\kappa(u_i,u_j)}\bigg)^\gamma\,\mathrm dx_j
    =\kappa(u_i,u_j)M_{\varphi,\gamma}.
\]
Iterating over the tree yields
\(
    M_{\varphi,\gamma}^k\prod_{\{i,j\}\in E(\mathcal T)}\kappa(u_i,u_j).
\)
Since $d_i$ is the degree of vertex $i$, Condition~(\condPK) gives
\(
    \prod_{\{i,j\}\in E(\mathcal T)}\kappa(u_i,u_j)
    \le c_\kappa^k\prod_{i=0}^k u_i^{d_i},
\)
which proves the claim.
\end{proof}

\def\tPP{\widetilde{\mathcal{P}}}
\def\Pn{\mathcal P _n}
\def\Tn{W_n}
\section{Proof for polynomial test functions}
\label{sec:coupling_discretization}

In this section, we prove the CLT for polynomial test functions stated in Theorem \ref{thm:polynomials}. For functionals of the Poisson process alone, optimal convergence rates are typically derived from the general results of \cite{mehler} stated in terms of finite difference operators when inserting one or more fixed points.  However, when working with functionals of an RCM, one needs to keep track of how edges are added to the extra points. We therefore use 
normal approximation bounds for edge-marked Poisson processes developed in \cite{nestmann} which are analogues of \cite{mehler} but with the point-wise difference operators replaced by  difference operators for edge marked Poisson processes introduced in Section \ref{ss:diff}.

 Section \ref{subsec:nestmann} records the Poincar\'e inequality, Wasserstein, and Kolmogorov bounds for normal approximations  \cite{nestmann}. Section \ref{subsec:new2} then proves the required moment and spatial integrability estimates for the first- and second-order insertion differences and combines them with the variance asymptotics from Section \ref{subsec:variance_comparison} to obtain both the quantitative normal approximation and the asserted central limit theorem.

\subsection{The CLT from \cite{nestmann}}\label{subsec:nestmann}

Recall the construction of the RCM and the difference operators from  Section~\ref{ss:diff}. For a functional $F_n= H(G(\xi_n))$, we have the following Poincaré inequality from \cite{nestmann}.
\begin{theorem}[Thm. 5.2 in \cite{nestmann}]\label{thm:poincare}
Suppose $\E F^2_n<\infty$. Then,
\(\Var F_n \le \int_{\check W_n} \E [(\Delta_{\check x_1; \check \xx_1} F_n)^2] \lambda(\d \check x_1).\)
\end{theorem}

Moreover, we have the following CLT:
\begin{theorem}[Thm. 6.1 in \cite{nestmann}] \label{thm:nestmann}
Suppose $\E F_n=0$, $\Var F_n=1$ and $\E F_n^4< \infty$. Let $\gamma_{1},\ldots,\gamma_6$ be as in \eqref{eq:gammas} below. If $\gamma_5,\gamma_6<\infty$  and the condition \eqref{eq:finite} below holds, then
$$d_W(F_n,\mathcal N(0,1)) \le \gamma_1 + \gamma_2 + \gamma_3$$
and 
$$d_K(F_n,\mathcal  N(0,1))\le \gamma_1+\gamma_2+\gamma_3 + \gamma_4 + \gamma_5 + \gamma_6.$$
\end{theorem}

The quantities $\gamma_i$ have the usual form from \cite{mehler}, but they must be interpreted through the  difference operators for edge-marked Poisson processes. 
Precisely,
\begin{align}\nonumber
&\gamma_1 = 2\bigg(\int_{\check W_n^3} \E[(\Delta_{\check x_1;\check \xx_3}F_n)^2(\Delta_{\check x_2;\check \xx_3}F_n )^2]^{1/2}\E[(\Delta^2_{\check x_1,\check x_3;\check \xx_3}F_n)^2(\Delta^2_{\check x_2,\check x_3;\check\xx_3} F_n )^2]^{1/2} \lambda^3(\d (\check x_1,\check x_2,\check x_3))\bigg)^{1/2}\\ \nonumber
&\gamma_2=\bigg(\int_{\check W_n^3} \E[(\Delta^2_{\check x_1,\check x_3;\check \xx_3}F_n)^2(\Delta^2_{\check x_2,\check x_3; \check \xx_3} F_n )^2] \lambda^3(\d (\check x_1,\check x_2,\check x_3))\bigg)^{1/2}\\ \nonumber
&\gamma_3=\int_{\check W_n} \E[|\Delta_{\check x_1;\check \xx_1} F_n|^3] \lambda(\d \check x_1) \\ \nonumber
\end{align}
and
\begin{align}
\nonumber
&\gamma_4=\frac{1}{2} \E [F^4_n] ^{1/4} \int_{\check W_n} \E[(\Delta_{\check x_1;\check \xx_1} F_n)^4]^{3/4} \lambda(\d \check x_1) \\ \nonumber
&\gamma_5=\bigg(\int_{\check W_n} \E[(\Delta_{\check x_1;\check \xx_1} F_n)^4] \lambda(\d \check x_1)\bigg)^{1/2}\\
&\gamma_6=\bigg(\int_{\check W_n^2}( 6\E[(\Delta_{\check x_1;\check \xx_2}F_n)^4]^{1/2}\E[(\Delta^2_{\check x_1,\check x_2;\check \xx_2}F_n )^4]^{1/2} + 3\E[(\Delta^2_{\check x_1,\check x_2;\check \xx_2}F_n)^4]) \lambda^2(\d (\check x_1,\check x_2))\bigg)^{1/2}.\label{eq:gammas}
\end{align}
Moreover, the following two finiteness conditions are required:
\begin{align} \nonumber
& \int_{\check W_n^3} \E[(\Delta_{\check x_1;\check \xx_3}F_n)^4]^{1/4} \E [(\Delta_{\check x_2;\check \xx_3}F_n )^4]^{1/4}\E[(\Delta^2_{\check x_1,\check x_3;\check \xx_3}F_n)^4]^{1/4} \E[(\Delta^2_{\check x_2,\check x_3;\check \xx_3} F_n )^4]^{1/4} \lambda^3(\d (\check x_1,\check x_2,\check x_3))<\infty,\\
&\int_{\check W_n^3} \E[(\Delta^2_{\check x_1,\check x_3;\check \xx_3}F_n)^4]^{1/2} \E[(\Delta^2_{\check x_2,\check x_3;\check \xx_3} F_n )^4]^{1/2} \lambda^3(\d (\check x_1,\check x_2,\check x_3))<\infty. \label{eq:finite}
\end{align} 
To bound $\E[F_n^4]$ in $\gamma_4$ we may use:
\bel[Lemma 6.3 in \cite{nestmann}]\label{lem:F4} If $\E [F^4_n]<\infty$, $\E F_n =0 $ and $\Var (F_n)=1$, then
$$ \E [F_n^4] \leq \max \bigg\{ 256\bigg(\int_{\check W_n}  \E[(\Delta_{\check x_1; \check \xx_1}F_n)^4]^{1/2}\lambda (\d \check x_1)\bigg)^2, 4\int_{\check W_n} \E [(\Delta_{\check x_1; \check \xx_1} F_n)^4] \lambda (\d \check x_1) + 2 \bigg\}. $$
\enl

In the following, set
\begin{align}\label{deff}
 F_n:= F(\xi_n) := L(f;G(\xi_n)),
\end{align}
and whenever $\operatorname{Var}(F_n)>0$,  
define $$\widehat F_n :=\frac{F(\xi_n) - \E[F(\xi_n)]}{\sqrt{\textup{Var}(F(\xi_n))}}.
$$ 
Then we have
$F_n\in L^4$ and hence $\widehat F_n\in L^4$ for every fixed $n$. Indeed, letting $A_n$ denote the adjacency matrix of
$G(\xi_n)$, the bounds $\|A_n\|_{\mathrm{op}}\le N_n$
and $|f(t)|\le C(1+|t|^m)$ give
$|F_n|\le CN_n(1+N_n^m)$, and all moments of the
Poisson variable $N_n$ are finite. 
We may thus apply Theorem \ref{thm:nestmann} to $\widehat F_n$ and get
\begin{align} \label{clt}
	d_W\big(\widehat F_n, \mc N(0,1) \big) \le \textup{Var}(F_n)^{-1}{\g_1} + \textup{Var}(F_n )^{-1}{\g_2} + \textup{Var}(F_n )^{-3/2}\g_3.
\end{align} 
and
\begin{align} \label{cltK}
	d_K\big(\widehat F_n, \mc N(0,1) \big) \le \textup{Var}(F_n)^{-1}{\g_1} + \textup{Var}(F_n )^{-1}{\g_2} + \textup{Var}(F_n )^{-3/2}\g_3 + \hat\g_4
    + \textup{Var}(F_n )^{-1}\g_5 
    + \textup{Var}(F_n )^{-1}\g_6
\end{align} 
where the $\gamma_i$ are computed for $F_n$ and $\hat\gamma_4$ is computed for $\hat F_n$.
The powers of the variance in these bounds follow directly
from the homogeneity of the difference operators.
Indeed, writing $s_n:=\sqrt{\operatorname{Var}(F_n)}>0$, 
\[
\Delta_{\check x_1;\check \xx_1}\widehat F_n
=s_n^{-1}\Delta_{\check x_1;\check \xx_1}F_n,
\qquad
\Delta^2_{\check x_1,\check x_2;\check \xx_2}\widehat F_n
=s_n^{-1}\Delta^2_{\check x_1,\check x_2;\check \xx_2}F_n.
\]
Inspection of \eqref{eq:gammas} therefore gives $\gamma_i(\widehat F_n)
=s_n^{-2}\gamma_i(F_n) $ for $i\in\{1,2,5,6\}$, 
and $\gamma_3(\widehat F_n)
=s_n^{-3}\gamma_3(F_n).$ 
For $\gamma_4$, the factor
$\mathbb E[\widehat F_n^4]^{1/4}$ must also be controlled,
so we keep the notation $\widehat\gamma_4$ and estimate it
separately below.
It remains to control the moments of the first- and second-order difference operators of $F_n$.

%
%
\subsection{Quantitative CLT}
\label{subsec:new2}

To prove Theorem \ref{thm:polynomials}, we first establish moment bounds for the first- and second-order add-one cost operators $\Delta_{\check x_1;\check \xx_1}F_n$ and $\Delta^2_{\check x_1,\check x_2;\check \xx_2}F_n$. 
Set
\[
\overline M_\varphi:=1\vee M_{\varphi,1},\qquad
\Psi_{\ell,m}(x,y):=\mathbb E\bigg[\varphi\bigg(\frac{\|x-y\|^d}{(\ell m)^d\kappa(U,V)}\bigg)U^{\ell m-1}V^{\ell m-1}\bigg],
\]
where $U,V$ are independent random variables with distribution $\mu$. 

\begin{lemma}[Uniform fourth-moment bounds for difference operators]
\label{lem:uniform_moments}
Let $f$ be a polynomial of degree $m\ge1$. Assume that \eqref{cond:PK} holds, that \(M_{\varphi,1/4}<\infty\), and that \(\mathbb E[U^{4m-1}]<\infty\). 
There exists a constant $C<\infty$, independent
of $n$, such that the following bounds hold
for every $n\ge1$.
\been
\im[(i)]  \textbf{First Difference Operator:} For every $\check \xx_1:= \check x=(x,u)\in\check W_n$,
    \begin{align}
        \mathbb{E}\big[|\Delta_{(x,u); \check \xx_1} F_n|^4\big] \le C u^{4m-1}\overline M_\varphi^{4m-1}.
    \end{align}

\im[(ii)]  \textbf{Second Difference Operator:}  For every pair of distinct spatial points
$x_1,x_2\in W_n$,
    \begin{align}
        \iint\mathbb{E}\big[|\Delta_{(x_1,u_1),(x_2,u_2); \check \xx_2}^2 F_n|^4 \big] \mu(\d u_1) \mu(\d u_2) \le C \Psi_{4,m}(x_1,x_2).
    \end{align}
\enen
\end{lemma}

The spatial integrations in the normal-approximation bounds use the following scale-mixture estimate.

\begin{lemma}[Scale-mixture integrability]
\label{lem:scale_mixture}
Let $\gamma\in(0,1]$, $\ell,m\ge1$, and set $a=\ell m-1$. Assume that \eqref{cond:PK} holds, that $M_{\varphi,\gamma}<\infty$, and that $\mathbb E[U^b]<\infty$ for some $b>a$ satisfying $\gamma(b-a)>1$. Then,
\[
\int_{\mathbb R^d}\Psi_{\ell,m}(0,z)^\gamma\,\mathrm dz<\infty.
\]
\end{lemma}

\begin{proof}
  Since $\kappa(U,V)\le c_\kappa UV$, decompose according to the events $2^k\le UV<2^{k+1}$, $k\ge0$. Monotonicity of $\varphi$ gives
\[
\Psi_{\ell,m}(0,z)
\le
\mathbb E[(UV)^b]\sum_{k=0}^{\infty} 2^{-k(b-a)}
\varphi\bigg(\frac{\|z\|^d}{2c_\kappa (\ell m)^d2^k}\bigg).
\] 
For $0<\gamma\le1$ and nonnegative numbers $a_k$, we have $\left(\sum_k a_k\right)^\gamma
\le
\sum_k a_k^\gamma.$ 
Applying this inequality to the preceding bound for
$\Psi_{\ell,m}(0,z)$, we have
\[
\begin{aligned}
\int_{\mathbb R^d}\Psi_{\ell,m}(0,z)^\gamma\,dz
&\le
\mathbb E[(UV)^b]^\gamma
\sum_{k=0}^{\infty}2^{-k\gamma(b-a)} 
\int_{\mathbb R^d}
\varphi\!\left(
 \frac{\|z\|^d}{2c_\kappa (\ell m)^d2^k}
\right)^\gamma dz.
\end{aligned}
\]
For the  integral, setting $s_k:=2c_\kappa (\ell m)^d2^k,$ 
\[
\begin{aligned}
\int_{\mathbb R^d}
\varphi\!\left(\frac{\|z\|^d}{s_k}\right)^\gamma dz
&=
s_k\int_{\mathbb R^d}\varphi(\|y\|^d)^\gamma\,dy =
2c_\kappa (\ell m)^d2^kM_{\varphi,\gamma}.
\end{aligned}
\]
It follows that
\[
\begin{aligned}
\int_{\mathbb R^d}\Psi_{\ell,m}(0,z)^\gamma\,dz
&\le
2c_\kappa (\ell m)^d
\mathbb E[U^b]^{2\gamma}
M_{\varphi,\gamma}
\sum_{k=0}^{\infty}
2^{-k(\gamma(b-a)-1)}.
\end{aligned}
\]
The geometric series converges because $\gamma(b-a)>1$.
This proves the lemma.
\end{proof}

Given Lemma \ref{lem:uniform_moments}, we conclude the proof of Theorem \ref{thm:polynomials}.
\bep[Proof of Theorem \ref{thm:polynomials}] 
Let $f$ be a polynomial of degree $m$, and let $F_n$ be defined in
\eqref{deff}. 
By Lemma~\ref{lem:uniform_moments},
\begin{equation}\label{eq:A-bound}
A_x(u):=\E\big[|\Delta_{(x,u);\check \xx_1}F_n|^4\big]\le C u^{4m-1},
\end{equation}
and for \(B_{x_1,x_2}(u_1,u_2):=
\E\big[|\Delta^2_{(x_1,u_1),(x_2,u_2);\check \xx_2}F_n|^4\big]\), we have
\begin{equation}\label{eq:B-bound}
\iint B_{x_1,x_2}(u_1,u_2)\,\mu(du_1)\mu(du_2)
\le C\Psi_{4,m}(x_1,x_2).
\end{equation}
Assume $m\ge2$.
We first bound $\gamma_1$.
For fixed inserted positions and marks, Cauchy--Schwarz  gives 
\[
\mathbb E\!\left[
 (\Delta_{\check x_1;\check \xx_3}F_n)^2
 (\Delta_{\check x_2;\check \xx_3}F_n)^2
\right]^{1/2}
\le
A_{x_1}(u_1)^{1/4}A_{x_2}(u_2)^{1/4},
\]
and similarly
\[
\begin{aligned}
&\mathbb E\!\left[
 (\Delta^2_{\check x_1,\check x_3;\check \xx_3}F_n)^2
 (\Delta^2_{\check x_2,\check x_3;\check \xx_3}F_n)^2
\right]^{1/2} \le
B_{x_1,x_3}(u_1,u_3)^{1/4}
B_{x_2,x_3}(u_2,u_3)^{1/4}.
\end{aligned}
\] 
Substituting these inequalities into the definition of $\gamma_1$
gives
\[
\begin{aligned}
\frac{\gamma_1^2}{4}
\le
\int_{W_n^3}\int
&A_{x_1}(u_1)^{1/4}A_{x_2}(u_2)^{1/4} B_{x_1,x_3}(u_1,u_3)^{1/4}
B_{x_2,x_3}(u_2,u_3)^{1/4}
\,\mu ^{3} (\d(u_1,u_2,u_3))  \ dx_1dx_2dx_3.
\end{aligned}
\]
For fixed $x_1,x_2,x_3$, apply H\"older's inequality with
four exponents equal to $4$ on the common mark space:
\[
\begin{aligned}
&\int
 A_{x_1}(u_1)^{1/4}A_{x_2}(u_2)^{1/4}
 B_{x_1,x_3}(u_1,u_3)^{1/4}
 B_{x_2,x_3}(u_2,u_3)^{1/4}\,d\mu ^{3}\\
&\quad\le
\left(\int A_{x_1}(u_1)\,d\mu ^{3}\right)^{1/4}
\left(\int A_{x_2}(u_2)\,d\mu ^{3}\right)^{1/4} 
\left(\int B_{x_1,x_3}(u_1,u_3)\,d\mu ^{3}\right)^{1/4}
\left(\int B_{x_2,x_3}(u_2,u_3)\,d\mu ^{3}\right)^{1/4}.
\end{aligned}
\]
Integration over
a mark absent from an integrand contributes a factor $1$.
Thus, by \eqref{eq:A-bound}
\[
\int A_{x_i}(u_i)\, d\mu ^{3}
=
\int A_{x_i}(u_i)\,\mu(du_i)
\le C\mathbb E[U^{4m-1}]<\infty,
\]
whereas \eqref{eq:B-bound} gives
\[
\int B_{x_i,x_3}(u_i,u_3)\, d\mu ^{3}
\le C\Psi_{4,m}(x_i,x_3),
\qquad i=1,2.
\]
Therefore,
\[
\gamma_1^2
\le
C\int_{W_n^3}
\Psi_{4,m}(x_1,x_3)^{1/4}
\Psi_{4,m}(x_2,x_3)^{1/4}
\,dx_1dx_2dx_3.
\] 
For $m\ge2$, recalling the condition $\mathbb E[U^{6m}]<\infty$, Lemma~\ref{lem:scale_mixture} applies with
$\ell=4$, $\gamma=1/4$, and $b=6m$, since
\(
(6m-(4m-1))={2m+1}>4.
\)  
Hence by translation invariance of $\Psi_{4,m}$,
\[
\begin{aligned}
\gamma_1^2 \le C\int_{W_n}
\left(\int_{W_n}\Psi_{4,m}(x_1,x_3)^{1/4}\,dx_1\right)
\left(\int_{W_n}\Psi_{4,m}(x_2,x_3)^{1/4}\,dx_2\right)
dx_3 \le C|W_n|.
\end{aligned}
\] 
For $\gamma_2$, we similarly obtain
\[
\begin{aligned}
\gamma_2^2
&\le
\int_{W_n^3}
\int
B_{x_1,x_3}(u_1,u_3)^{1/2}
B_{x_2,x_3}(u_2,u_3)^{1/2}
\,\mu^3(d(u_1,u_2,u_3))
\,d(x_1,x_2,x_3)\\
&\le
C\int_{W_n^3}
\Psi_{4,m}(x_1,x_3)^{1/2}
\Psi_{4,m}(x_2,x_3)^{1/2}
\,d(x_1,x_2,x_3),
\end{aligned}
\]
where the second inequality follows from Cauchy--Schwarz in the
three mark variables. Since
$M_{\varphi,1/2}\le M_{\varphi,1/4}$, Lemma
\ref{lem:scale_mixture} with $\gamma=1/2$ and $b=6m$ yields
\(
\gamma_2^2\in O( |W_n|).
\)
 The same estimates also verify the
two finiteness conditions in \eqref{eq:finite}.
Next, by \eqref{eq:A-bound}, we have
\(
\gamma_3
\le
\int_{\check W_n}A_x(u)^{3/4}\,\lambda(d(x,u))
\in O( |W_n|),
\)
since $\E[U^{3(4m-1)/4}]<\infty$.

We now consider the further terms needed for the Kolmogorov
distance. Here, \(\int_{\check W_n}A_x(u)\,\lambda(d(x,u)) \in O(|W_n|)\) and
\(
\int_{\check W_n}A_x(u)^{1/2}\,\lambda(d(x,u))
\in O( |W_n|).
\)
If $\sigma_f^2>0$, then by Theorem~\ref{thm:variance_limit},
$\Var(F_n)\asymp |W_n|$.  The two integrals in Lemma~\ref{lem:F4} satisfy
\[
\int_{\check W_n}
\mathbb E[(\Delta_{\check x;(\check x)}\widehat F_n)^4]^{1/2}
\,\lambda(d\check x)
\le C\frac{|W_n|}{\Var(F_n)},\qquad \int_{\check W_n}
\mathbb E[(\Delta_{\check x;(\check x)}\widehat F_n)^4]
\,\lambda(d\check x)
\le C\frac{|W_n|}{\Var(F_n)^2}.
\] Lemma~\ref{lem:F4} therefore gives
\(
\E[\widehat F_n^4]\le C.
\)
It follows that
\[
\begin{aligned}
\widehat\gamma_4
&=
\frac12\mathbb E[\widehat F_n^4]^{1/4}
\Var(F_n)^{-3/2}
\int_{\check W_n}A_x(u)^{3/4}\,\lambda(d\check x) \le C\Var(F_n)^{-3/2}|W_n|
\le C|W_n|^{-1/2}.
\end{aligned}
\]
Moreover,
\(
\gamma_5^2
=
\int_{\check W_n}A_x(u)\,\lambda(d(x,u))
\in O( |W_n|).
\)

Finally, consider $\gamma_6$. By Cauchy--Schwarz in the two mark
variables,
\[
\begin{aligned}
&\iint
A_{x_1}(u_1)^{1/2}
B_{x_1,x_2}(u_1,u_2)^{1/2}
\,\mu(du_1)\mu(du_2)\\
&\qquad\le
\bigg(\int A_{x_1}(u_1)\,\mu(du_1)\bigg)^{1/2}
\bigg(\iint B_{x_1,x_2}(u_1,u_2)
\,\mu(du_1)\mu(du_2)\bigg)^{1/2} \le C\Psi_{4,m}(x_1,x_2)^{1/2},
\end{aligned}
\]
whereas directly from \eqref{eq:B-bound}, we get that
\( \iint B_{x_1,x_2}(u_1,u_2)
\,\mu(du_1)\mu(du_2)
\le C\Psi_{4,m}(x_1,x_2).
\)
It follows that
\[
\gamma_6^2
\le
C\int_{W_n^2}
\left(
\Psi_{4,m}(x_1,x_2)^{1/2}
+\Psi_{4,m}(x_1,x_2)
\right)\,d(x_1,x_2)
\in O( |W_n|),
\]
by Lemma~\ref{lem:scale_mixture}  with
$\gamma=1/2$ and $\gamma=1$, respectively.

If $m\in\{0,1\}$, then $F_n=a_0N_n$, because $\operatorname{Tr}(A_n)=0$.
It follows that $\Delta_{\check x;(\check x)} F_n=a_0$ and $\Delta_{\check x_1,\check x_2;\check \xx_2}^2 F_n=0$, from which we obtain $A_x(u)\le a_0^4$ and $B_{x_1,x_2}(u_1,u_2)=0$. This can be inserted in the estimates above to obtain (at most) the same order of the $\gamma_i$.

Assume first that $\sigma_f^2>0$. Theorem~\ref{thm:variance_limit} gives $\Var(F_n)\ge\sigma_f^2|W_n|/2$ for all sufficiently large $n$. Inserting the preceding estimates into \eqref{clt} yields
\(
d_W(\widehat F_n,\mathcal N(0,1))\le C|W_n|^{-1/2}.
\)
For the Kolmogorov distance,  the preceding estimates give 
$\widehat\gamma_4\le C|W_n|^{-1/2}$ and $\gamma_5+\gamma_6\le C|W_n|^{1/2}$, so their normalized contributions in \eqref{cltK} are also at most $C|W_n|^{-1/2}$, yielding $d_K(\widehat F_n,\mathcal N(0,1))\le C|W_n|^{-1/2}$.
Moreover,
\(
(F_n-\mathbb E[F_n])/{\sqrt{|W_n|}}
=
\widehat F_n\sqrt{\Var(F_n)/|W_n|},
\)
and Slutsky's theorem together with the variance limit proves convergence to $\mathcal N(0,\sigma_f^2)$.
If $\sigma_f^2=0$, then we have
\(
\mathbb E\bigg[\bigg|(F_n-\mathbb E[F_n])/{\sqrt{|W_n|}}\bigg|^2\bigg]
=
{\Var(F_n)}/{|W_n|}\longrightarrow0.
\)
Thus, the volume-normalized statistic converges in distribution to $\mathcal N(0,0)$, which completes the proof.
\enp

It remains to show Lemma \ref{lem:uniform_moments}. 
Observe that, by \eqref{eq:sumscore}, $F_n$ can be expressed as a sum of scores
\begin{align} 
\label{sumofscore}
  L(f;G(\xi_n)) =
\sum_{\check X_i\in\check P_n}
g_f(\check X_i,G(\xi_n)),
\end{align}
where the score function $g_f$ is given for a finite graph $G$ with marked vertex set $\check Q\subseteq \R^d\times [1,\infty)$  and a  marked vertex
$\check x_1=(x_1,u_1)\in\check Q$ by 
\begin{align}\label{score}
  g_f(\check x_1, G) &:= a_0 + \sum_{p = 2}^m  \sum_{\substack{(\check x_2,\dots,\check x_p) \in  (\check Q\setminus \{\check x_1\} )^{(p-1)}}} \sum_{q = p}^ma_q {\f1{(p-1)!}} \sum_{\substack{\pi \in \Sigma_{q,p}^\circ \\ \pi(1)=1} }\prod_{j =1}^{q} A(\check x_{\pi(j)}, \check x_{\pi(j + 1)} ; G).
\end{align}
A non-zero term can be interpreted as a closed walk of length $q$ starting and ending at $\check{x}_1$ visiting exactly $p$ distinct points. 
Therefore, there are only $(p-1)!$ enumerations of the remaining
vertices.
This explains the coefficient $1/(p-1)!$ in  the above expression.
Summing over all possible starting points gives that
 \eqref{sumofscore} is equivalent to \eqref{eq:sumscore}. 

To derive Lemma \ref{lem:uniform_moments}, we first establish corresponding moment bounds for the scores $g_f$ in Lemma \ref{lem:mom2}. We begin by introducing  notations. 
Fix an ordered tuple
\[
\check{\xx}_r=(\check x_1,\ldots,\check x_r),
\qquad
\check x_i=(x_i,u_i)\in\mathbb R^d\times[1,\infty),
\]
with  distinct spatial coordinates.
For $I\subseteq\{1,\ldots,r\}$, write
\[
G_{n,I,\check{\xx}_r}:=G_{n,I}(\xi,\check{\xx}_r),
\qquad
G_{n,\check{\xx}_r}
:=G_{n,\{1,\ldots,r\},\check{\xx}_r} 
\]
(see Section~\ref{ss:diff} for the definition of $G_{n,I}$).

We now define deletion differences for score functions.
For a finite graph $G$ containing distinct marked vertices
$\check z$ and $\check x$, define the first-order difference 
\[
\partial_{\check x}g_f(\check z,G)
:=
g_f(\check z,G)
-
g_f(\check z,G\setminus\{\check x\}).
\]
Deleting a marked vertex means deleting its spatial
vertex and all incident edges, while retaining every
other vertex, mark, and edge.
Similarly for  distinct marked vertices $\check z,\check x,\check y$
of $G$, define the second-order difference 
\[
\begin{aligned}
\partial^2_{\check x,\check y}g_f(\check z,G)
:={} 
g_f(\check z,G)
-g_f(\check z,G\setminus\{\check x\}) -g_f(\check z,G\setminus\{\check y\})
+g_f(\check z,G\setminus\{\check x,\check y\}).
\end{aligned}
\]
For example, 
\[
\begin{aligned}
\partial_{\check x_i}
g_f(\check z,G_{n,\check \xx_r})
={} 
g_f(\check z,G_{n,\{1,\ldots,r\},\check \xx_r}) -
g_f(\check z,G_{n,\{1,\ldots,r\}\setminus\{i\},
                    \check \xx_r}).
\end{aligned}
\]

\bel[Moment bounds]
\label{lem:mom2}
Let {$\ell, m \in \bN$} and let $f$ be a polynomial of degree $m$. Assume that \eqref{cond:PK} holds, that $M_{\varphi,1}<\infty$, and that $\mathbb E[U^{\ell m -1}]<\infty$.
Let $\check x_1 = (x_1, u_1) ,\check x_2 = (x_2, u_2) \in \check\T_n$ with $x_1 \neq x_2$.
Then, there exists a constant $C >0$ {(depending on $\ell$ and $f$, but not on $n,\check x_1,\check x_2$)} such that, for every $n \in \bN$, the following holds:
\been
\im[(i)] 
$\E[|g_f(\check x_1, G_{n, \check \xx_1} )|^{\ell}]\le C u_1^{\ell m-1}\overline M_\varphi^{\ell m-1}.
$
\im[(ii)] 
\(
\E \Big[ \big\vert \sum_{\check X_i \in \check \PP_n} \partial_{\check x_1} g_f(\check X_i,G_{n,  \check \xx_1})\big\vert^\ell\Big]\le C u_1^{\ell m-1}\overline M_\varphi^{\ell  m-1}.
\)
\item[(iii)] 
\(\iint
\E\left[
|\partial_{\check x_1}
g_f(\check x_2 ,G_{n,  \check \xx_2 })|^\ell
\right]
\mu(du_1)\mu(du_2)
\le
C\overline M_\varphi^{\ell m-2}
\Psi_{\ell,m}(x_1,x_2).\)
\im[(iv)] 
\(
\iint
\E\Big[
\Big|
\sum_{\check X_i\in\check\PP_n}
\partial^2_{\check x_1, \check x_2}
g_f(\check X_i,G_{n,\check\xx_2})
\Big|^\ell 
\Big]
\mu(du_1)\mu(du_2)
\le
C\overline M_\varphi^{\ell m-2}\Psi_{\ell,m}(x_1,x_2).
\)
\enen 
\enl

\bep
Write $f(x)=\sum_{q=0}^m a_qx^q$ and let $g_q:=g_{x^q}$ denote the score associated with $x^q$. The score and its difference operators are linear in $f$, and for every finite family $(Y_q)_{q=0}^m$, we have
\(
\Big|\sum_{q=0}^m a_qY_q\Big|^\ell\le C_{f,\ell}\sum_{q=0}^m|Y_q|^\ell.
\)
It is therefore enough to establish the estimates for a monomial. The constant term has zero score differences, and the linear term contributes nothing because the adjacency matrices have zero diagonal. We henceforth take $f(x)=x^m$ with $m\ge 2$, for which the score in \eqref{score} is given by
 \begin{align} \label{polyscore}
   g_f(\check x_1, G ) &:= \sum_{p = 2}^m  \sum_{\substack{(\check x_2,\dots,\check x_p) \in (\check Q\setminus \{\check x_1\} )^{(p-1)}}} {\f1{(p-1)!}} \sum_{\substack{\pi \in \Sigma_{m,p}^\circ \\ \pi(1)=1} }\prod_{j =1}^{m} A(\check x_{\pi(j)},\check x_{\pi(j + 1)} ; G),
 \end{align}
We prove the four parts of the lemma separately.\\[-2ex]

\emph{Part (i).} Since we have reduced to $f(x)=x^m$, the score $g_f(\check x_1,G)$ counts rooted closed walks of length $m$ starting and ending at $\check x_1$. Consequently, its $\ell$-th power counts ordered $\ell$-tuples of such walks. Concatenating the $\ell$ walks produces a closed walk of length $\ell m$ starting and ending at $\check x_1$. This map is injective, since the original walks are recovered by cutting the concatenated walk after every $m$ steps. Hence, writing $\check X_{-1}:=\check x_1$, $i_1:=-1$ and $\pi(\ell m +1) := \pi(1)$, 
\begin{align} \label{501}
  \big|g_f(\check x_1, G_{n, \check \xx_1})\big|^\ell \le \sum_{p = 2}^{\ell m}  \sum_{\substack{(\check X_{i_2},\dots,\check X_{i_p}) \in  \check \PP_{n} ^{(p-1)}}} {\f1{(p-1)!}} \sum_{\substack{\pi \in \Sigma_{\ell m,p}^\circ \\ \pi(1)=1} }\prod_{j =1}^{\ell m} A(\check X_{i_{\pi(j)}}, \check X_{i_{\pi(j + 1)}} ; G_{n, \check \xx_1})
\end{align}

For $\pi\in\Sigma^\circ_{\ell m,p}$, consider the
connected simple undirected graph formed by the
distinct edges traversed by the corresponding
closed walk.  For each $\pi$, choose its spanning tree $\mathcal{T}_\pi$. Since every entry of the adjacency matrix is either 0 or 1,  
\begin{align} \label{502}
  \prod_{j=1}^{\ell m} A\bigl(\check X_{i_{\pi(j)}}, \check X_{i_{\pi(j+1)}};\, G_{n, \check \xx_1} \bigr)
 \le 
\prod_{\{j,k\}\in E(\mc T_\pi)} A\bigl(\check X_{i_j}, \check X_{i_k};\,  G_{n, \check \xx_1} \bigr),
\end{align}
where \(E(\mc T_\pi)\) denotes the collection of (undirected) edges in $\mc T_\pi$ and \(\{j,k\}\) represents the (undirected) edge between vertices \(j\) and \(k\).
 Thus, using \eqref{501} and \eqref{502}, 
\begin{align} \label{505}
  & |g_f(\check x_1, G_{n, \check \xx_1})|^\ell\le C_{\ell,m}\sum_{p = 2}^{\ell m} \sum_{(\check X_{i_2},\dots,
  \check X_{i_p}) \in \check \PP_{n}^{(p-1)}} {\f1{(p-1)!}} \sum_{\substack{\pi \in \Sigma_{\ell m,p}^\circ \\ \pi(1)=1}}\prod_{ \{j,k\}\in E(\mc T_\pi)} A(\check X_{i_j},\check X_{i_k}; G_{n, \check \xx_1}).
\end{align}
For a fixed $\pi \in \Sigma_{\ell m,p}^\circ$, the multivariate Mecke formula gives, writing $\check x_k = (x_k, u_k)$ for $k\in \{2,\dots,p\},$
\begin{align} \label{506}
&\E\Big[\sum_{(\check X_{i_2},\dots,\check X_{i_p}) \in \check \PP_{n}^{(p-1)}} \prod_{\{j,k\}\in E(\mc T_\pi)} A(\check X_{i_j},\check X_{i_k} ; G_{n, \check \xx_1}) \Big ] \nonumber \\
&= \int_{\T_n^{p-1}} \int_{[1,\infty)^{p-1}} \E\Big[ \prod_{\{j,k\}\in E(\mc T_\pi) } A(\check x_{j} ,\check x_{k} ;  G_{n, \check \xx_p}) \Big] \prod_{k=2}^p \d x_k  \prod_{k=2}^p \mu(\d u_k).
\end{align} 
Recalling the connection rule, this  is equal to
\begin{align} \label{508}  \int_{\T_n^{p-1}} \int_{[1,\infty)^{p-1}} \prod_{ \substack{\{j,k\}\in E(\mc T_\pi) \\ k<j}} \vp\Big(\frac{\|x_j-x_k\|^d}{\kappa(u_j,u_k)}\Big)\prod_{k=2}^p \d x_k  \prod_{k=2}^p \mu(\d u_k).
\end{align}
The argument now proceeds by integrating over leaf vertices of the tree $\mc T_\pi$. Lemma~\ref{lem:weighted_tree_estimate}, applied conditionally on the marks, bounds the spatial integral by
\(
\overline M_\varphi^{p-1}c_\kappa^{p-1}u_1^{d_1}\prod_{i=2}^p u_i^{d_i},
\)
where $d_i$ is the degree of vertex $i$ in $\mc T_\pi$. Since $d_i\le p-1\le\ell m-1$, direct integration of the mark powers gives
\[
\mathbb E\Big[\sum_{(\check X_{i_2},\dots,\check X_{i_p})\in\check\PP_n^{(p-1)}}\prod_{\{j,k\}\in E(\mc T_\pi)}A(\check X_{i_j},\check X_{i_k};G_{n,\check\xx_1})\Big]
\le C\overline M_\varphi^{p-1}u_1^{\ell m-1},
\]
provided that $\mathbb E[U^{\ell m-1}]<\infty$. 
Returning to \eqref{505}, take expectations and sum this bound
over the finitely many values of $p$ and $\pi$:
\[
\begin{aligned}
\mathbb E|g_f(\check x_1,G_{n, \check \xx_1})|^\ell
&\le
Cu_1^{\ell m-1}
\sum_{p=2}^{\ell m}
\frac{|\Sigma_{\ell m,p}^\circ|}{(p-1)!}
\overline M_\varphi^{\,p-1} \le
C' u_1^{\ell m-1}
\overline M_\varphi^{\,\ell m-1}.
\end{aligned}
\]
This proves Part (i). 
\medskip

%
%
\emph{Part (ii).}  
Fix $\check X_i\in\check P_n$ and set $i_1:=i$.
Recall that $\check X_{-1}:=\check x_1$ is the
inserted vertex.
Write $G_n:=G(\xi_n)$ and
$\check \PP_{n,\check x_1}:=\check \PP_n\cup\{\check x_1\}$.
By the definition of the score deletion difference,
\[
\partial_{\check x_1}g_f(\check X_{i},G_{n, \check \xx_1})
=
g_f(\check X_{i},G_{n, \check \xx_1})
-
g_f(\check X_{i},G_n).
\]
Using \eqref{polyscore} with $\pi(m+1)=\pi(1)$, we have
\[
g_f(\check X_{i},G_{n})= \sum_{p = 2}^m  \sum_{\substack{(\check X_{i_2},\dots,\check X_{i_p}) \in  ( \check  \PP_{n} \setminus \{\check X_{i}\})^{(p-1)}}} {\f1{(p-1)!}} \sum_{\substack{\pi \in \Sigma_{m,p}^\circ \\ \pi(1)=1} }\prod_{j =1}^{m} A(\check X_{i_{\pi(j)}}, \check X_{i_{\pi(j + 1)}}; G_{n}) 
\]
and  
\[
g_f(\check X_{i}, G_{n, \check \xx_1})=
\sum_{p = 2}^m  \sum_{\substack{(\check X_{i_2},\dots,\check X_{i_p}) \in  ( \check \PP_{n, \check x_1}  \setminus \{\check X_{i}\}) ^{(p-1)}}} {\f1{(p-1)!}} \sum_{\substack{\pi \in \Sigma_{m,p}^\circ \\ \pi(1)=1} }\prod_{j =1}^{m} A(\check X_{i_{\pi(j)}}, \check X_{i_{\pi(j + 1)}} ;  G_{n, \check \xx_1} ).
\]
The induced subgraph of $G_{n, \check \xx_1}$ on the old
vertices is exactly $G_n$, with the same edge indicators.
Consequently, every term whose auxiliary tuple $(\check X_{i_2},\ldots,\check X_{i_p})$ avoids
$\check x_1$ cancels between these two expansions.
The remaining terms are precisely those whose auxiliary
tuple contains $\check x_1$.
They are nonnegative and count closed walks based at
$\check X_{i}$ that visit the inserted vertex.
Therefore,  
\begin{align*}
    \partial_{\check x_1} g_f(\check X_{i},G_{n, \check \xx_1} ) = \sum_{p = 2}^m  \sum_{\substack{(\check X_{i_2},\dots,\check X_{i_{p}}) \in  (\check \PP_{n,\check x_1}\setminus \{\check X_{i}\})^{(p-1)}\\{\check x_1 \in \{\check X_{i_2},\dots,\check X_{i_{p}}\}}}} {\f1{(p-1)!}} \sum_{\substack{\pi \in \Sigma_{m,p}^\circ \\ \pi(1)=1} }\prod_{j =1}^{m} A(\check X_{i_{\pi(j)}}, \check X_{i_{\pi(j + 1)}} ;  G_{n, \check \xx_1} ).
\end{align*}
Thus, 
\begin{align*}
&\Big(\sum_{\check X_i \in \check \PP_n} \partial_{\check x_1} g_f(\check X_i,G_{n, \check \xx_1} )\Big)^\ell \\
&\qquad = \prod_{h=1}^\ell \Big( \sum_{p_h = 2}^{m}  \sum_{\substack{(\check X_{i_{h,1}},\dots,\check X_{i_{h,p_h}}) \in  \check \PP_{n,\check x_1}^{(p_h)}\\{\check x_1 \in \{\check X_{i_{h,2}},\dots,\check X_{i_{h,p_h}}\}}}} {\f1{(p_h-1)!}} \sum_{\substack{\pi_h \in \Sigma^\circ_{m,p_h} \\ \pi_h(1)=1} }\prod_{j =1}^{m} A(\check X_{i_{h,\pi_h(j)}}, \check X_{i_{h,\pi_h(j + 1)}} ;  G_{n, \check \xx_1}) \Big).
\end{align*}
Indeed,
summing over $\check X_i\in\check\PP_n$ amounts to
including the starting vertex as the first coordinate
$\check X_{i_{h,1}}$ of each ordered $p_h$-tuple. 
Thus, the  summation over $\check X_i$ is
 incorporated into the summation over
$\check\PP_{n,\check x_1}^{(p_h)}$.

For each $h$, there is a unique position
$s_h\in\{2,\ldots,p_h\}$ occupied by $\check x_1$.
Separating the sum according to this position rewrites
the preceding right-hand side as
\begin{equation}
\begin{aligned}\label{eq:Dxbound}
\prod_{h=1}^{\ell}\Bigg[
&\sum_{p_h=2}^{m}\sum_{s_h=2}^{p_h}
\frac{1}{(p_h-1)!}
\sum_{\substack{\pi_h\in\Sigma^\circ_{m,p_h}\\\pi_h(1)=1}}
\sum_{(\check X_{i_{h,k}})_{k\in[p_h]\setminus\{s_h\}}
      \in\check \PP_n^{(p_h-1)}}
\prod_{j=1}^{m}
A\bigl(
 \widehat X_{h,\pi_h(j)},
 \widehat X_{h,\pi_h(j+1)};
 G_{n, \check \xx_1}
\bigr)
\Bigg],
\end{aligned}
\end{equation}
where $[p_h]:=\{1,\ldots,p_h\}$,
$\pi_h(m+1):=\pi_h(1)$, and
\[
\widehat X_{h,j}
:=
\begin{cases}
\check X_{i_{h,j}},&j\ne s_h,\\
\check X_{-1}=\check x_1,&j=s_h.
\end{cases}
\]
The tuple
 $(\check X_{i_{h,k}})_{k\in[p_h]\setminus\{s_h\}}$  is ordered by increasing
$k\in[p_h]\setminus\{s_h\}$.
For fixed  $p_h,s_h,\pi_h$,  write
\[
g_h\bigl((\check z_{h,k})_{k\in[p_h]\setminus\{s_h\}};G\bigr)
:=
\prod_{j=1}^{m}
A\bigl(
 \widehat z_{h,\pi_h(j)},
 \widehat z_{h,\pi_h(j+1)};G
\bigr),
\]
where $\widehat z_{h,j}:=\check z_{h,j}$ for $j\ne s_h$ and $\widehat z_{h,j}:=\check x_1$ for $j= s_h$.

Consider $\ell$ rows, where row $h$ corresponds to the
$h$th walk and contains $p_h-1$ Poisson vertices,
indexed by $(h,j)$ with $j\in[p_h]\setminus\{s_h\}$.
When expanding the product in \eqref{eq:Dxbound}, the Poisson
vertices within each row are distinct, but those from different
rows may coincide. We encode these coincidences by partitions
of  
\[
\mathcal J:=\{(h,j):1\le h\le\ell,\ 
j\in[p_h]\setminus\{s_h\}\}.
\]
Let $\Pi(p_1-1,\ldots,p_\ell-1)$ denote the set of partitions
$\sigma$ of $\mathcal J$ such that each block contains at most
one element from each row. Each block represents one distinct
Poisson vertex. Write $|\sigma|$ for the number of blocks,
order them by their lexicographically smallest elements,
and let  $b_\sigma(h,j)\in\{1,\ldots,|\sigma|\}$ be the index of
the block containing $(h,j)$.
For $\check\yy_{|\sigma|}=(\check y_1,\ldots,\check y_{|\sigma|})$,
set
\begin{equation}\label{eq:otimesh}
\left(\bigotimes_{h=1}^{\ell}g_h\right)_\sigma
(\check\yy_{|\sigma|};G)
:=
\prod_{h=1}^{\ell}
g_h\bigl(
 (\check y_{b_\sigma(h,j)})_{j\in[p_h]\setminus\{s_h\}};
 G
\bigr).
\end{equation}Also define $G_{\sigma,n}
:=
G_{\{1,\ldots,|\sigma|+1\}}
\bigl(
 \xi_n,
 (\check x_1,\check y_1,\ldots,\check y_{|\sigma|})
\bigr).$ 
By multivariate Mecke formula 
\cite[Theorem~3.1]{bogdan} ($\check y_1,\ldots,\check y_{|\sigma|}$ are  
integration variables),  
 the expectation of \eqref{eq:Dxbound} is given by 
\[
\begin{aligned}
&
\sum_{\substack{
(p_h,s_h,\pi_h)_{h=1}^{\ell}:\\
2\le p_h\le m,\;2\le s_h\le p_h,\\
\pi_h\in\Sigma^\circ_{m,p_h},\;
\pi_h(1)=1\ \text{for all }h
}}
\left(\prod_{h=1}^{\ell}\frac{1}{(p_h-1)!}\right)
\sum_{\sigma\in\Pi(p_1-1,\ldots,p_\ell-1)} 
\int_{\check W_n^{|\sigma|}}
\mathbb E\!\left[
 \left(\bigotimes_{h=1}^{\ell}g_h\right)_\sigma
 (\check\yy_{|\sigma|};G_{\sigma,n})
\right]
\prod_{a=1}^{|\sigma|}\lambda(d\check y_a)
\end{aligned}
\] 
(recall $\lambda(d\check y)=dy\,\mu(dv)$).
Using the block ordering defined above, write $\check w_1:=\check x_1$ and $\check w_{a+1}:=\check y_a$ ($1\le a\le |\sigma|$). 
Let $H_\sigma$ be the (undirected) graph on
$\{1,\ldots,|\sigma|+1\}$ whose edges are the distinct edges
traversed by the rows after the identifications
prescribed by $\sigma$.
Choose its spanning tree $T_\sigma$.  Then,
\begin{align*}
    \left(\bigotimes_{h=1}^{\ell}g_h\right)_\sigma
 (\check\yy_{|\sigma|};G_{\sigma,n}) \le   \prod_{\{j,k\} \in E(\mc T_\sigma)} A(\check w_j,\check w_k;G_{\sigma, n}).
\end{align*} 
The tree $\mathcal T_\sigma$ has $|\sigma|$ edges.
Let $d_j$ be the degree of vertex $j$ in $\mathcal T_\sigma$. By conditionally averaging the distinct edge indicators and applying Lemma~\ref{lem:weighted_tree_estimate}, we obtain 
\[
\begin{aligned}
&\int_{\check W_n^{|\sigma|}}
\mathbb E\!\left[
 \left(\bigotimes_{h=1}^{\ell}g_h\right)_\sigma
 (\check{\yy}_{|\sigma|};G_{\sigma,n})
\right]
\prod_{a=1}^{|\sigma|}\lambda(d\check y_a)\le
c_\kappa^{|\sigma|} M_\varphi^{|\sigma|}
u_1^{d_1}
\prod_{j=2}^{|\sigma|+1}\mathbb E[U^{d_j}].
\end{aligned}
\] 
The root position $x_1$ and root mark $u_1$ remain fixed.
Since $|\sigma|\le\sum_{h=1}^{\ell}(p_h-1)
\le\ell(m-1)\le\ell m-1$ 
and $d_j\le |\sigma|$, the moment assumption  implies the bound above is at most $C_{\ell,m}
u_1^{\ell m-1}\overline M_\varphi^{\,\ell m-1}$.   
Summing these bounds in the preceding diagram formula
proves Part (ii).

\medskip

\medskip
\noindent\emph{Part (iii).}
Set $L:=\ell m$.
For $f(x)=x^m$, the quantity $\partial_{\check x_1}g_f(\check x_2,G_{n, \check \xx_2})$ 
counts length-$m$ closed walks based at $x_2$ that visit
$x_1$.
It is nonnegative, and its $\ell$th power counts ordered
collections of $\ell$ such walks.  Concatenating these walks at $x_2$ produces a closed walk
of length $L$ visiting both $x_1$ and $x_2$.
Let $p$ denote the number of distinct vertices visited by
this concatenated walk. Then $2\le p\le L$.
Label these vertices by $\{1,\ldots,p\}$, assigning labels
$1$ and $2$ to $x_1$ and $x_2$, respectively.
Let $H$ be the simple (undirected)  graph on $\{1,\ldots,p\}$
whose edges are  traversed by
the walk.

Apply the multivariate Mecke formula to the $p-2$ vertices
labeled $3,\ldots,p$, which belong to the original Poisson
process, while keeping the two distinguished vertices
$\check x_1$ and $\check x_2$ fixed.  
Then average the independent edges.
After also integrating the two inserted marks, the
contribution of this pattern is
\[
\int_{[1,\infty)^p}\int_{W_n^{p-2}}
\prod_{\{i,j\}\in E(H)}
\varphi\!\left(
 \frac{\|z_i-z_j\|^d}{\kappa(u_i,u_j)}
\right)
\,dz_3\cdots dz_p\prod_{v=1}^p\mu(du_v),
\]
where $z_1=x_1$ and $z_2=x_2$ are fixed. 
There are finitely many patterns, depending only on
$\ell,m$, so it suffices to bound each such integral.
Fix a spanning tree $\mathcal T$ of the  graph $H$.
Write
\[
p_{\{i,j\}}
:=
\varphi\!\left(
 \frac{\|z_i-z_j\|^d}{\kappa(u_i,u_j)}
\right).
\]
Since these factors lie in $[0,1]$, the product over
$E(H)$ is bounded by the product over $E(\mathcal  T)$.

Let $P_{\mathcal  T}(x_1,x_2)$ be the unique path joining $x_1$ and $x_2$. 
For every spatial configuration, the triangle inequality
therefore guarantees at least one path edge of length
at least $\|x_1-x_2\|/L$. Thus
\[
\prod_{h\in E(\mathcal  T)}p_h
\le
\sum_{e=\{a,b\}\in P_{\mathcal  T}(x_1,x_2)}
\mathbf1_{\{\|z_a-z_b\|\ge\|x_1-x_2\|/L\}}
\prod_{h\in E(\mathcal  T)}p_h.
\]
For each fixed summand, monotonicity of $\varphi$ gives $\mathbf1_{\{\|z_a-z_b\|\ge\|x_1-x_2\|/L\}}p_e
\le
\varphi (
 \frac{\|x_1-x_2\|^d}{L^d\kappa(u_a,u_b)}
 ).$ 
Consequently,
\[
\prod_{h\in E(\mathcal  T)}p_h
\le
\sum_{e=\{a,b\}\in P_{\mathcal  T}(x_1,x_2)}
\varphi\!\left(
 \frac{\|x_1-x_2\|^d}{L^d\kappa(u_a,u_b)}
\right)
\prod_{h\in E(\mathcal  T)\setminus\{e\}}p_h.
\]
Fix an edge $e=\{a,b\}$ in this sum.
Its deletion leaves two trees, each containing exactly
one fixed root.
After enlarging the spatial domain to
$(\mathbb R^d)^{p-2}$, integrate all nonroot vertices
successively from the leaves of these two trees.
Lemma~\ref{lem:weighted_tree_estimate}  gives the bound
\[
c_\kappa^{p-2}M_\varphi^{p-2}
\varphi\!\left(
 \frac{\|x_1-x_2\|^d}{L^d\kappa(u_a,u_b)}
\right)
\prod_{v=1}^p u_v^{d'_v},
\]
where $d'_v$ is the degree in the forest $\mathcal T\setminus\{e\}$. 
The endpoint labels $a,b$ and the forest are now fixed
within this summand.
Independence of  marks gives
\[
\begin{aligned}
&\int_{[1,\infty)^p}
\varphi\!\left(
 \frac{\|x_1-x_2\|^d}{L^d\kappa(u_a,u_b)}
\right)
\prod_{v=1}^p u_v^{d'_v}
\prod_{v=1}^p\mu(du_v) =
\left(
 \prod_{v\notin\{a,b\}}\mathbb E[U^{d'_v}]
\right)
\mathbb E\!\left[
 \varphi\!\left(
  \frac{\|x_1-x_2\|^d}{L^d\kappa(U,V)}
 \right)
 U^{d'_a}V^{d'_b}
\right],
\end{aligned}
\] where $U,V$ are independent random variables with distribution $\mu$. 
As $d'_v\le p-1\le L-1$,
the first factor is   bounded by a finite
constant.
For the second factor, as $U,V\ge1$,
\[
\begin{aligned}
&\mathbb E\!\left[
 \varphi\!\left(
  \frac{\|x_1-x_2\|^d}{L^d\kappa(U,V)}
 \right)
 U^{d'_a}V^{d'_b}
\right] \le
\mathbb E\!\left[
 \varphi\!\left(
  \frac{\|x_1-x_2\|^d}{L^d\kappa(U,V)}
 \right)
 U^{L-1}V^{L-1}
\right]
=
\Psi_{\ell,m}(x_1,x_2).
\end{aligned}
\] 
Finally,
$M_\varphi^{p-2}\le\overline M_\varphi^{\,L-2}$.
Summing over the finitely many path edges and 
patterns  concludes the proof.

\medskip

\medskip
\noindent\emph{Part (iv).}
Continue with $f(x)=x^m$, $m\ge2$, and set
$L:=\ell m$ and $G:=G_{n,\check\xx_2}$.
Write
\[
S_m:=
\sum_{\check z\in\check\PP_n}
\partial^2_{\check x_1,\check x_2}g_f(\check z,G).
\]
For a vertex $\check z$, the second score
difference counts length-$m$ closed walks based
at $z$ that visit both inserted vertices.
Indeed, every walk missing at least one inserted
vertex cancels in the alternating deletion
difference.
Thus $S_m\ge0$, and $S_m^\ell $ counts ordered
$\ell $-tuples of such walks.
Unlike in Part (iii), their starting vertices
are old Poisson vertices and need not coincide.

Group these trajectories according to the equality
patterns of their $\ell m$ indexed walk positions.
Each pattern distinguishes the two blocks occupied
by $x_1,x_2$, both of which meet every walk;
each starting position belongs to a nonroot block.
Following the notation in Part (iii), let $p$
be the number of distinct visited vertices,
label $x_1,x_2$ by $1,2$, and label the remaining
vertices in order of first appearance across
the ordered walks.
Let $H$ be the simple undirected graph formed
by the distinct traversed edges.
Then $p\le L=\ell m$, and $H$ is connected.

For each fixed pattern, average the edge randomness
and apply the multivariate Mecke formula
as in Part (iii).
After also integrating the two inserted marks,
its contribution is the same two-root graph
integral considered there, with $z_1=x_1$,
$z_2=x_2$, and $p-2$ Poisson vertices.  
Applying the aforementioned argument  to bound each
pattern's contribution by
\[
C_m\overline M_\varphi^{\,\ell m-2}
\Psi_{\ell,m}(x_1,x_2).
\]
The required moment 
$\mathbb E[U^{\ell m-1}]<\infty$ is included in
the assumptions.
Since the number of equality patterns depends only
on $\ell,m$, summing their contributions concludes
the proof.

\enp

Finally, we proceed with the proof of Lemma \ref{lem:uniform_moments}, inspired by the argument of \cite[Lemma 5.5]{mal_stab}. 
%
%
\bep[Proof of Lemma \ref{lem:uniform_moments}] We prove the two parts separately.
\medskip

\noindent \emph{Part (i).}
Recalling the definition of the difference operator, we get as in \cite[Lemma 5.2]{mal_stab} that
\begin{align}
	\label{lem:451}
	\Delta_{\check x_1;\check \xx_1}F_n = g_f(\check x_1,G_{n, \check \xx_1}) + \sum_{\check X_i \in\check \PP_n} \partial_{\check x_1} g_f( \check X_i, G_{n,\check \xx_1}).
\end{align}
Thus,
$
	\E \big[|\Delta_{\check x_1;\check \xx_1}F_n|^4 \big]
 \le 16 \E\big[ g_f(\check x_1,G_{n, \check \xx_1} )^4\big] + 16 \E\Big[\Big( \sum_{\check X_i \in\check \PP_n} \partial_{\check x_1} g_f(\check X_i, G_{n,\check \xx_1})\Big)^4\Big].
 $
By Lemma \ref{lem:mom2}(i), the first term is controlled as
$
\E\big[ g_f(\check x_1,G_{n, \check \xx_1} )^4\big] \le C u_1^{4m-1} \overline M_\varphi^{4m-1}.
$
By Lemma \ref{lem:mom2}(ii), the second term is bounded by
\[
\mathbb E\Big\vert 
 \sum_{\check X_i\in\check P_n}
 \partial_{\check x_1}
 g_f(\check X_i,G_{n,\check \xx_1})
\Big\vert ^4
\le
C u_1^{4m-1}\overline M_\varphi^{\,4m-1}.
\] 
This proves Part (i).

\emph{Part (ii).} 
Recalling the definition of second-order difference operator, we have
$$
\Delta_{\check x_1,\check x_2;\check \xx_2}^2 F_n=\partial_{\check x_1}g_f(\check x_2,G_{n, \check \xx_2} ) + \partial_{\check x_2}g_f (\check x_1, G_{n, \check \xx_2} ) + \sum_{\check X_i \in\check \PP_n} \partial_{\check x_1,\check x_2}^2 g_f ( \check X_i, G_{n, \check \xx_2} )
$$
(see \cite[Lemma 5.2]{mal_stab} for details),
implying that
\begin{align*}
 &\int_{[1,\infty)^2} 
 \mathbb E\left|
  \Delta^2_{(x_1,u_1),(x_2,u_2);\check \xx_2}F_n
 \right|^4
 \,\mu (\d u_1)\mu (\d u_2) \nonumber \\ & \le C\int_{[1,\infty)^2}\Big(\E |\partial_{\check x_1}g_f(\check x_2,G_{n, \check \xx_2}) |^4 
 + \E |\partial_{\check x_2}g_f(\check x_1,G_{n, \check \xx_2}) |^4 +\E\Big\vert  \sum_{\check X_i \in\check \PP_n} \partial_{\check x_1, \check x_2 }^2 g_f(\check X_i,G_{n,\check  \xx_2})\Big\vert ^4 \Big) \mu(\d u_1)\mu(\d u_2). 
\end{align*}
 By Lemma \ref{lem:mom2} (iii), the first and second terms are both bounded by
\(
C \overline M_\varphi^{4m-2} \Psi_{4,m}(x_1,x_2),
\)
where $C>0$ is a constant. 
For the last term, we obtain by Lemma \ref{lem:mom2} (iv),
$$
\int_{[1,\infty)^2} \E\Big[\Big(\sum_{\check X_i \in\check \PP_n} \partial_{\check x_1,\check x_2}^2 g_f(\check X_i,G_{n,\check  \xx_2})\Big)^4\Big] \mu(\d u_1) \mu(\d u_2) \le C\overline M_\varphi^{4m-2}\Psi_{4,m}(x_1,x_2).
$$
This gives the assertion.
 
\enp

%
%
\section{Proof for general test functions}
\label{sec:general}

We now extend the CLT for polynomial test functions in Theorem \ref{thm:polynomials} to the  broader class of test functions $\mathcal{H}_{c,0}^2$. This relies on an abstract approximation principle stated in Proposition \ref{aux prop} below that allows us to pass from polynomials to general test functions. The argument is essentially the same as in
\cite[Section~5.1]{hirsch-nam-otto-rgg} and is not specific to the present model. 

In Section \ref{ss:key} we state the abstract approximation principle and show how Theorem \ref{thm:general_functions} follows through two key intermediate results. These key results are proved in Section \ref{ss:key_proof} with the main difference for the weighted RCM being that the variance estimate (Proposition \ref{main prop}) is derived from the RCM version of the  Poincaré inequality stated in Theorem \ref{thm:poincare} involving the  RCM difference operator
\(\Delta_{\check x}\) rather than the point-wise difference operator for Poisson input. Otherwise, the proof follows the Fourier-Duhamel argument of \cite[Thm. 5.2]{hirsch-nam-otto-rgg}. The only 
model-specific input  is an exponential trace estimate stated separately   (Lemma \ref{lem:exponential_trace_bounds}) obtained by orienting
edges according to the vertex marks.

\subsection{Key ingredients and proof of Theorem \ref{thm:general_functions}}
\label{ss:key}

We first state the following abstract approximation principle used to pass from polynomial
test functions to general test functions. All constants in the estimates below are independent of \(n\),
unless explicitly stated otherwise.

\begin{proposition}[Abstract approximation principle]
\label{aux prop}
Let \(\mc L\) be a real vector space of functions equipped with a norm
\(\|\cdot\|_{\mc L}\), and let \(Z_n:\mc L\to L^2(\Omega)\) be linear maps.
Assume that there is a dense subspace \(\widetilde{\mc L}\subseteq\mc L\) and a
quadratic form \(V:\widetilde{\mc L}\to\mathbb R_{\ge0}\) such that
\(Z_n(\varphi)\xrightarrow{d}\mc N(0,V(\varphi))\) for every
\(\varphi\in\widetilde{\mc L}\). Assume also that, for some \(C>0\),
\begin{equation}
\label{eq:abstract_var}
    \sup_{n\ge1}\mathbb E[Z_n(\varphi)^2]
    \le
    C\|\varphi\|_{\mc L}^2
    \qquad\text{for all }\varphi\in\mc L .
\end{equation}
Then, \(V\) extends uniquely to a continuous quadratic form on \(\mc L\), and
\(Z_n(\varphi)\xrightarrow{d}\mc N(0,V(\varphi))\) for every
\(\varphi\in\mc L\).
\end{proposition}

This is similar to \cite[Prop. 5.1]{hirsch-nam-otto-rgg} adapted from \cite{cltsparse}. We give a proof for the convenience of the reader.

\begin{proof}
For \(\varphi,\psi\in\widetilde{\mc L}\), the second moment bound and the assumed
convergence on \(\widetilde{\mc L}\) give
\[
    V(\varphi-\psi)
    \le
    \liminf_{n\to\infty}\mathbb E[Z_n(\varphi-\psi)^2]
    \le
    C\|\varphi-\psi\|_{\mc L}^2 .
\]
Hence, \(V\) is continuous on \(\widetilde{\mc L}\), and it extends uniquely by
continuity to \(\mc L\).

Fix \(\varphi\in\mc L\), and choose \(\varphi_k\in\widetilde{\mc L}\) with
\(\|\varphi_k-\varphi\|_{\mc L}\to0\). The same variance bound gives
\[
    \sup_{n\ge1}
    \mathbb E\big[|Z_n(\varphi)-Z_n(\varphi_k)|^2\big]
    \le
    C\|\varphi-\varphi_k\|_{\mc L}^2 .
\]
For each fixed \(k\), \(Z_n(\varphi_k)\xrightarrow{d} \mc N(0,V(\varphi_k))\), and
\(V(\varphi_k)\to V(\varphi)\). 
Let \(Y_k\sim\mc N(0,V(\varphi_k))\) and
\(Y\sim\mc N(0,V(\varphi))\). If \(h\) is bounded by \(1\) and
\(1\)-Lipschitz, then
\[
\begin{split}
    \big|\mathbb E h(Z_n(\varphi))-\mathbb E h(Y)\big|
    \le{}&
    \mathbb E\big[|Z_n(\varphi)-Z_n(\varphi_k)|\big]  
    +
    \big|\mathbb E h(Z_n(\varphi_k))-\mathbb E h(Y_k)\big|
    +
    \big|\mathbb E h(Y_k)-\mathbb E h(Y)\big|.
\end{split}
\]
Taking \(\limsup_{n\to\infty}\) and then letting \(k\to\infty\) shows weak convergence.
\end{proof}

To use the abstract approximation to show Theorem  \ref{thm:general_functions}, we need two input results.
 The first input is a  lemma that allows us to use Proposition \ref{aux prop} to pass from  polynomials to $\mathcal{H}_{c,0}^2$.

\begin{lemma}[Polynomial density in \(\mathcal H_{c,0}^2\)]
\label{lem:polynomial_density_Lc}
Let \(c>0\). Then, polynomials \(p\) satisfying \(p(0)=0\) are dense in
\(\mathcal H_{c,0}^2\), that is, for every \(f\in\mathcal H_{c,0}^2\), there are polynomials
\(p_k\) with \(p_k(0)=0\) such that
\[
    \|f-p_k\|_{\mathcal H_c^2}\to0 .
\]
\end{lemma}

The second input  controls the variance of exponentially twisted statistics and allows us to verify \eqref{eq:abstract_var}.

\begin{proposition}[Variance estimate for exponentially twisted statistics]
\label{main prop}
Assume that \eqref{cond:EM} and \eqref{cond:PK} hold and that \(M_{\varphi,1}<\infty\).
For every \(c\ne0\), there exists \(C>0\), depending only on \(c\) and the
model parameters, such that for every \(f\in H^2(\mathbb R)\)   satisfying
\(f(0)=0\),
\[
    \Var\big(\Tr[f(A_n)e^{cA_n}]\big)
    \le
    C|W_n|\big(\|f\|_2^2+\|f''\|_2^2\big).
\]
\end{proposition}

 The two preceding inputs now give the CLT for the general test-functions.
\begin{proof}[Proof of Theorem~\ref{thm:general_functions}]
We apply Proposition~\ref{aux prop} with \(\mathcal L=\mathcal H_{c,0}^2\). By Lemma~\ref{lem:polynomial_density_Lc}, we may take $\tilde{\mathcal{L}}$ to be the space of polynomials \(p\) satisfying
\(p(0)=0\). For
\(f\in\mathcal H_{c,0}^2\), define
\[
    Z_n(f)
    :=
    |W_n|^{-1/2}
    \big(\Tr[f(A_n)]-\mathbb E\Tr[f(A_n)]\big).
\]
 By Theorem~\ref{thm:polynomials},
\(Z_n(p)\) converges to a centered Gaussian random variable for every such
polynomial \(p\). The limiting variance $\sigma_p^2$ is a quadratic form in \(p\) since it is the limit of quadratic forms. 

It remains to verify \eqref{eq:abstract_var}. Let \(\psi(x):=f(x)\sech(cx)\). Then,
\(\psi(0)=0\) and \(f(x)=\psi(x)\cosh(cx)\), and hence 
\(
    \Tr[f(A_n)]
    =
    (\Tr[\psi(A_n)e^{cA_n}]
    +
    \Tr[\psi(A_n)e^{-cA_n}])/2 .
\)
Since $\psi\in H^2(\R)$, we may apply
Proposition~\ref{main prop} to $\psi$ with constants \(c\) and \(-c\). Together with the inequality
\(\Var((X+Y)/2)\le(\Var X+\Var Y)/2\), this gives
\[
    |W_n|^{-1}\Var(\Tr[f(A_n)])
=|W_n|^{-1}\Var( (\Tr[\psi(A_n)e^{cA_n}]
    +
   \Tr[\psi(A_n)e^{-cA_n}])/2 )\le
    C\big(\|\psi\|_2^2+\|\psi''\|_2^2\big).
\]
If \(h(x):=\sech(cx)\), then \(|h'|\le |c|h\) and \(|h''|\le c^2h\). Since 
\(\psi=fh\),
\[
    \|\psi\|_2+\|\psi''\|_2
    \le
    C_c\big(\|fh\|_2+\|f'h\|_2+\|f''h\|_2\big)
    \le
    C_c\|f\|_{\mathcal H_c^2}.
\]
Hence,
\(
    \sup_{n\ge1}\mathbb E[Z_n(f)^2]
    \le
    C\|f\|_{\mathcal H_c^2}^2
    \)
    {for all } \(f\in\mathcal H_{c,0}^2.\)
Proposition~\ref{aux prop} therefore extends the polynomial CLT to every
\(f\in\mathcal H_{c,0}^2\) and Theorem \ref{thm:general_functions} follows.
\end{proof}

\subsection{Proof of key ingredients}
\label{ss:key_proof}
We now prove Lemma \ref{lem:polynomial_density_Lc}.

\begin{proof}[Proof of Lemma \ref{lem:polynomial_density_Lc}]

Fix $c>0$ and set $w(x):=\operatorname{sech}^2(cx).$ 
Then $e^{-2c|x|}\le w(x)\le4e^{-2c|x|}.$

We first prove that the set of polynomials is dense in $L^2(w\,dx)$.
Suppose that $g\in L^2(w\,dx)$ is orthogonal to every
polynomial.
For $|\operatorname{Re}z|<c$, define
\[
F(z):=\int_{\mathbb R}g(x)w(x)e^{zx}\,dx.
\] 
Then $F$ is holomorphic on
$\{z:|\operatorname{Re}z|<c\}$, and orthogonality gives $F^{(k)}(0)
=
\int_{\mathbb R}x^kg(x)w(x)\,dx=0$ for all $k \ge 0.$
Hence
$F\equiv0$.
In particular, the Fourier transform of $gw$ vanishes, implying $gw=0$ a.e.
As $w>0$, we have $g=0$ a.e.
The orthogonal complement of the polynomials is
 trivial, proving density in $L^2(w\,dx)$.

Now, define $(Jg)(x):=\int_0^x g(t)\,dt.
$ 
Then $e^{-cx}(Jg)(x)
=
\int_0^x e^{-c(x-t)}e^{-ct}g(t)\,dt.$ 
The convolution kernel
$e^{-cs}\mathbf1_{\{s\ge0\}}$ has $L^1$ norm $1/c$. By
Young's convolution inequality,
\[
\int_0^\infty e^{-2cx}|Jg(x)|^2\,dx
\le
c^{-2}\int_0^\infty e^{-2cx}|g(x)|^2\,dx.
\]
By symmetry we have the corresponding inequality on
the negative domain. Hence, 
\begin{align} \label{intin}
\|Jg\|_{L^2(w)}
\le\frac2c\|g\|_{L^2(w)}.
\end{align}

Let $f\in H_c^2$.
Choose polynomials $q_k$ such that
$q_k\to f''$ in $L^2(w\,dx)$, and set
\[
p_k(x)
:=
f(0)+f'(0)x+\int_0^x(x-t)q_k(t)\,dt.
\]
Each $p_k$ is a polynomial, with
$p_k(0)=f(0)$, $p_k'(0)=f'(0)$, and $p_k''=q_k$.
Thus
\[
(f-p_k)''=f''-q_k,\qquad
(f-p_k)'=J(f''-q_k),\qquad
f-p_k=J^2(f''-q_k).
\]
Hence using \eqref{intin},
\[
\|f-p_k\|_{H_c^2}^2
\le
\left(1+\frac4{c^2}+\frac{16}{c^4}\right)
\|f''-q_k\|_{L^2(w)}^2
\longrightarrow0.
\]
If $f(0)=0$, then $p_k(0)=0$ for every $k$.
This concludes the proof.

\end{proof}

We next prove exponential trace bounds that are needed for the
Fourier--Duhamel argument in the proof of Proposition \ref{main prop}. 
 Recall that $A_n$ denotes the adjacency matrix of the weighted RCM $G(\xi_n)$. Let $D_n$ denote the diagonal matrix containing the degrees  of $G(\xi_n)$. For
$\check x\in\check W_n$, let $A_n^{\check x}$ be the adjacency matrix of $G_{n,\check x}:=G_{\{1\} }(\xi_n,(\check x))$. We write $\mathbf e_{\check x}$ for the standard basis vector corresponding to the $\check x$-coordinate and set
\[
\mathbf a_{\check x}:=A_n^{\check x}\mathbf e_{\check x}.
\]

\begin{lemma}[Trace bounds]
\label{lem:exponential_trace_bounds}
Assume that \eqref{cond:EM}  and \eqref{cond:PK} hold and that \(M_{\varphi,1}<\infty\). For every \(R>0\),
\begin{align}
\label{eq:uniform-exponential-trace}
    \sup_{n\ge1}\sup_{0\le r\le R}
    |W_n|^{-1}\mathbb E\big[\Tr(e^{rA_n})+\Tr(e^{-rA_n})\big]
    &<\infty,\\
\label{eq:uniform-degree-trace}
    \sup_{n\ge1}\sup_{0\le r\le R}
    |W_n|^{-1}\mathbb E\big[\Tr(D_nA_n^2e^{rA_n})\big]
    &<\infty.
\end{align}
Moreover,
\begin{align}
\label{eq:inserted-root-exponential-trace}
    \sup_{n\ge1}\sup_{0\le r\le R}
    |W_n|^{-1}
    \int_{\check W_n}
    \mathbb E\!\left[
        \left\langle e^{rA_n^{\check x}}\mathbf e_{\check x},
        \mathbf e_{\check x}\right\rangle
    \right]\lambda(\d\check x)
    &<\infty,\\
\label{eq:inserted-root-degree-trace}
    \sup_{n\ge1}\sup_{0\le r\le R}
    |W_n|^{-1}
    \int_{\check W_n}
    \mathbb E\!\left[
        \left\langle e^{rA_n^{\check x}}\mathbf a_{\check x},
        \mathbf a_{\check x}\right\rangle
        \left\langle\mathbf a_{\check x},
        \mathbf a_{\check x}\right\rangle
    \right]\lambda(\d\check x)
    &<\infty.
\end{align}
\end{lemma}

\begin{proof}
We first derive the  estimates  \eqref{eq:inserted-root-exponential-trace} and \eqref{eq:inserted-root-degree-trace}, given \eqref{eq:uniform-exponential-trace} and \eqref{eq:uniform-degree-trace}. By the Mecke formula \eqref{eq:rcm-mecke}, 
\[
\begin{aligned} \mathbb E\operatorname{Tr}(e^{rA_n}) =
\mathbb E\sum_{x\in \mathcal P_n}
\langle e^{rA_n}e_x,e_x\rangle  =
 \int_{\check W_n}
\mathbb E[
\langle e^{rA_n^{\check x}}e_{\check x},e_{\check x}\rangle]
\,\lambda(d\check x) 
.
\end{aligned}
\] 
Thus, \eqref{eq:inserted-root-exponential-trace} follows  from
\eqref{eq:uniform-exponential-trace}.
To prove \eqref{eq:inserted-root-degree-trace}, observe that
\(
\left\langle
\mathbf a_{\check x},\mathbf a_{\check x}
\right\rangle
=
d_{\check x}
\)
and, since
\(\mathbf a_{\check x}=A_n^{\check x}\mathbf e_{\check x}\), we have
\( 
\left\langle
e^{rA_n^{\check x}}\mathbf a_{\check x},
\mathbf a_{\check x}
\right\rangle
=
\left(
(A_n^{\check x})^2e^{rA_n^{\check x}}
\right)_{\check x\check x}.
\)
  Mecke formula therefore gives
\[
\begin{aligned}
&\int_{\check W_n}
\mathbb E\left[
\left\langle
e^{rA_n^{\check x}}\mathbf a_{\check x},
\mathbf a_{\check x}
\right\rangle
\left\langle
\mathbf a_{\check x},\mathbf a_{\check x}
\right\rangle
\right]\lambda(\d\check x) 
=
\mathbb E\sum_{x\in\mathcal P_n}
d_x\big(A_n^2e^{rA_n}\big)_{xx}
=
\mathbb E[\Tr(D_nA_n^2e^{rA_n})].
\end{aligned}
\]
The bound \eqref{eq:inserted-root-degree-trace} now follows from
\eqref{eq:uniform-degree-trace}.

~

We now verify \eqref{eq:uniform-exponential-trace} and  \eqref{eq:uniform-degree-trace}.
A direct way to bound the adjacency matrix is through total
degrees:
\(
|z^\top A_nz|
\le
\sum_{x\in\PP_n}d_xz_x^2
=
z^\top D_nz.
\)
The resulting exponential trace estimate would require
control of exponential moments of the total degree.
The degree of a vertex with
mark $u$ is Poisson with mean at most $C_0u$, where
$C_0:=c_\kappa M_{\varphi,1}\mathbb E U$ with $c_\kappa$ defined in \eqref{cond:PK}.
Consequently, 
for fixed $R>0$,
\[
\mathbb E[e^{Rd_x}\mid U_x=u]
\le
\exp\bigl(C_0(e^R-1)u\bigr).
\]
Averaging this bound would require
$\mathbb E[e^{C_0(e^R-1)U}]<\infty$.
Condition \eqref{cond:EM}  guarantees a finite  exponential moment only at some $\eta$, so this approach does not provide the desired estimate for arbitrary $R>0$.

We instead distinguish neighbors with smaller marks from
those with larger marks.
For a neighbor mark $v\ge u$, the factor $uv$, by the condition
\eqref{cond:PK},  is bounded by $v^2$.
Thus, a finite second mark moment already gives a uniform
bound for the conditional mean of neighbors with no smaller
mark.
The conditional mean of neighbors with smaller marks may
grow linearly in $u$.
Orienting edges toward larger marks separates these two
contributions.
The finite exponential moment assumption \eqref{cond:EM}  will be used for the
incoming contribution after its coefficient has been made
sufficiently small.

To specify the orientation also when marks are equal, attach
to each Poisson point an independent auxiliary variable
$T\sim\Unif[0,1]$, independently of everything else.
Write $(U_x,T_x)$ for the enlarged mark of vertex $x$, and orient
each existing edge from the smaller pair to the larger pair
in lexicographic order.
The variable $T$ is used only to break ties in $U$, whose
distribution may have atoms, including the case $U\equiv1$.
Thus, the adjacency matrix is not changed. 

For a vertex $x$, let $d_x^+$ and $d_x^-$ denote its outgoing
and incoming degrees, respectively, and let $D^+$ and $D^-$
be the corresponding diagonal matrices.
Then $D_n=D^++D^-$.
For $\tau>0$, set
\[
Q_\tau:=\tau D^++\tau^{-1}D^-.
\]
The asymmetric coefficients are chosen to match the different
moment bounds for the two degrees.
We will take $\tau$ large below.  

Now, for every oriented edge \(x\to y\), we have
\(
2|z_xz_y|
\le
\tau z_x^2+\tau^{-1}z_y^2.
\)
Summing this inequality over the oriented edges yields
\[
|z^\top A_nz|
\le
\sum_{x\to y}2|z_xz_y|
\le
\sum_{x\to y}\big(\tau z_x^2+\tau^{-1}z_y^2\big)
=
z^\top Q_\tau z.
\]
Thus, \(z^\top(\pm A_n)z\le z^\top Q_\tau z\) for every \(z\). 
Recall that if \(B, C\) are symmetric matrices of the same size and
\(
z^\top Bz\le z^\top Cz
\text{ for every }z,
\)
then
\(
\lambda_j(B)\le\lambda_j(C)
\)
for every \(j\);
see \cite[Corollary~III.2.3]{Bhatia1997}. Applying this result with
\(B=\pm A_n\) and \(C=Q_\tau\), we obtain
\(
\Tr(e^{rA_n})+\Tr(e^{-rA_n})
\le
2\Tr(e^{rQ_\tau})
=
2\sum_{x\in\mathcal P_n}
e^{r(\tau d_x^++\tau^{-1}d_x^-)}.
\)
Consequently,
\begin{equation}\label{eq:tr_bound}
\begin{aligned}
&|W_n|^{-1}
\mathbb E\big[\Tr(e^{rA_n})+\Tr(e^{-rA_n})\big]
\le
2|W_n|^{-1}
\mathbb E\sum_{x\in\mathcal P_n}
e^{r(\tau d_x^++\tau^{-1}d_x^-)}.
\end{aligned}
\end{equation}
It therefore remains to bound the averaged exponential moment of the
oriented degrees.

For $x\in W_n$, $u\ge1$, and $t\in[0,1]$, let
$\mathbb E^{x,u,t}$ denote expectation after inserting a vertex
at $x$ with enlarged mark $(u,t)$ into an independent marked
Poisson configuration in $W_n$ and sampling its incident edges
according to the connection rule.  
Independent Poisson thinning implies that $d_x^+$ and $d_x^-$ are 
independent Poisson variables under $\mathbb E^{x,u,t}$.
Write their means as $m_{n,+}(x,u,t)$ and $m_{n,-}(x,u,t)$.
We derive the bounds on these means explicitly.
Fix the inserted root $(x,u,t)$.
The remaining enlarged marked vertices form a Poisson process
on $W_n\times[1,\infty)\times[0,1]$ with intensity
$dy\,\mu(dv)\,ds$.
For a potential neighbor $(y,v,s)$, write
\[
p_{x,u}(y,v)
:=
\varphi\!\left(
 \frac{\|y-x\|^d}{\kappa(u,v)}
\right).
\] 
Conditional on the remaining marked configuration, the expected
contribution of $(y,v,s)$ to the incoming degree is $p_{x,u}(y,v)
\mathbf1_{\{(v,s)\prec_{\mathrm{lex}}(u,t)\}},$ 
and its expected contribution to the outgoing degree is $
p_{x,u}(y,v)
\mathbf1_{\{(u,t)\prec_{\mathrm{lex}}(v,s)\}}.$ 
Averaging these sums over the Poisson configuration  gives
\[
\begin{aligned}
m_{n,-}(x,u,t)
&=
\int_{W_n}\int_{[1,\infty)}\int_0^1
p_{x,u}(y,v)
\mathbf1_{\{(v,s)\prec_{\mathrm{lex}}(u,t)\}}
\,ds\,\mu(dv)\,dy,\\
m_{n,+}(x,u,t)
&=
\int_{W_n}\int_{[1,\infty)}\int_0^1
p_{x,u}(y,v)
\mathbf1_{\{(u,t)\prec_{\mathrm{lex}}(v,s)\}}
\,ds\,\mu(dv)\,dy.
\end{aligned}
\]  
The lexicographic order gives  
\[
\begin{aligned}
\int_0^1
\mathbf1_{\{(v,s)\prec_{\mathrm{lex}}(u,t)\}}\,ds \leq  \mathbf1_{\{v\le u\}}, \qquad 
\int_0^1
\mathbf1_{\{(u,t)\prec_{\mathrm{lex}}(v,s)\}}\,ds
 \leq \mathbf1_{\{v\ge u\}}.
\end{aligned}
\] 
For fixed $u,v\ge1$,  
\[
\begin{aligned}
\int_{W_n}p_{x,u}(y,v)\,dy
&\le
\int_{\mathbb R^d}
\varphi\!\left(
 \frac{\|y-x\|^d}{\kappa(u,v)}
\right)\,dy =
\kappa(u,v)
\int_{\mathbb R^d}\varphi(\|z\|^d)\,dz =
\kappa(u,v)M_{\varphi,1}.
\end{aligned}
\]  
Define
\[
\Lambda_-:=c_\kappa M_{\varphi,1}\mathbb E [U],
\qquad
\Lambda_+:=c_\kappa M_{\varphi,1}\mathbb E[U^2].
\]
These constants are finite under Condition \eqref{cond:EM}.
Using Condition \eqref{cond:PK}, we obtain
\begin{align*}
m_{n,-}(x,u,t)
\le
\int_{[1,u]}
 \left(\int_{W_n}p_{x,u}(y,v)\,dy\right)\mu(dv) &\le
M_{\varphi,1}\int_{[1,u]}\kappa(u,v)\,\mu(dv) \\
&\le
c_\kappa M_{\varphi,1}
u\int_{[1,u]}v\,\mu(dv) =
c_\kappa M_{\varphi,1}
u\mathbb E[U\mathbf1_{\{U\le u\}}]
\le\Lambda_-u.
\end{align*}
Similarly, for the outgoing degree,
\begin{align*}
m_{n,+}(x,u,t)
 \le
\int_{[u,\infty)}
 \left(\int_{W_n}p_{x,u}(y,v)\,dy\right)\mu(dv) 
&\le
M_{\varphi,1}\int_{[u,\infty)}\kappa(u,v)\,\mu(dv)\\
&\le
c_\kappa M_{\varphi,1}
u\mathbb E[U\mathbf1_{\{U\ge u\}}] 
 \le
c_\kappa M_{\varphi,1}
\mathbb E[U^2\mathbf1_{\{U\ge u\}}]
\le\Lambda_+.
\end{align*} 
The Poisson exponential-moment formula now gives, for
$\alpha,\beta\ge0$,
\[
\begin{aligned}
\mathbb E^{x,u,t}
 [e^{\alpha d_x^++\beta d_x^-}]
&=
\exp\!\left\{
 m_{n,+}(x,u,t)(e^\alpha-1)
 +m_{n,-}(x,u,t)(e^\beta-1)
\right\} \le
\exp\!\left\{
 \Lambda_+(e^\alpha-1)+\Lambda_-u(e^\beta-1)
\right\}.
\end{aligned}
\]
  $\{(X_i,U_i,T_i):1\le i\le N_n\}$ 
is Poisson with intensity $dx\,\mu(du)\,dt$.
By  Mecke formula,   
using the above bound,
\begin{equation}
 \begin{aligned}
 |W_n|^{-1}\mathbb E\sum_{x\in P_n}
e^{\alpha d_x^++\beta d_x^-} 
& =
|W_n|^{-1}
\int_{W_n}\int_{[1,\infty)}\int_0^1
\mathbb E^{x,u,t}[e^{\alpha d_x^++\beta d_x^-}]
\,dt\,\mu(du)\,dx\\
&\quad\le
|W_n|^{-1}
\int_{W_n}\int_{[1,\infty)}\int_0^1
e^{\Lambda_+(e^\alpha-1)+\Lambda_-u(e^\beta-1)}
\,dt\,\mu(du)\,dx\\
&\quad=
e^{\Lambda_+(e^\alpha-1)}
\mathbb E[e^{\Lambda_-U(e^\beta-1)}].
\end{aligned}
\end{equation} 
Choose $\tau_R\ge1$ sufficiently large that $\Lambda_-\bigl(e^{R/\tau_R}-1\bigr)\le\eta/2.$ 
For $0\le r\le R$, set $\alpha=r\tau_R$ and
$\beta=r/\tau_R$.
The right-hand side is then at most $\exp\!\left\{\Lambda_+\bigl(e^{R\tau_R}-1\bigr)\right\}
\mathbb E[e^{\eta U/2}]
 <\infty.$ 
 Inserting this into \eqref{eq:tr_bound} 
  proves \eqref{eq:uniform-exponential-trace}.

~

We next prove \eqref{eq:uniform-degree-trace}. For \(0\le r\le R\), set
\(
B_r:=A_n^2e^{rA_n}.
\)
The matrix \(B_r\) is positive semi-definite because its eigenvalues are
\(\lambda^2e^{r\lambda}\) where $\lambda$ are the eigenvalues of $A_n$. The Cauchy--Schwarz inequality (in Hilbert--Schmidt  norm) 
gives
\(\Tr(D_nB_r)
\le
\Tr(D_n^2)^{1/2}\Tr(B_r^2)^{1/2}.\)
Taking expectations and applying Cauchy--Schwarz, 
\begin{align} \label{44441}
\mathbb E[\Tr(D_nA_n^2e^{rA_n})]
\le
\mathbb E[\Tr(D_n^2)]^{1/2}
\mathbb E[\Tr(B_r^2)]^{1/2}.
\end{align} 
We bound each factor respectively.
For a root inserted at $x$ with mark $u$, its total degree
has Poisson mean
\[
\begin{aligned}
m_n(x,u)
&=
\int_{W_n}\int_{[1,\infty)}
\varphi\!\left(
 \frac{\|x-y\|^d}{\kappa(u,v)}
\right)\,\mu(dv)\,dy \le c_\kappa M_{\varphi,1}u\mathbb E U
=\Lambda_-u.
\end{aligned}
\]
This is a bound for the total mean itself. 
Thus, $\mathbb E^{x,u,t}[d_x^2]
\le\Lambda_-u+\Lambda_-^2u^2,$ 
uniformly also in the auxiliary root mark $t$.
Applying   Mecke formula,
\[
\begin{aligned}
\mathbb E\operatorname{Tr}(D_n^2)
&=
\int_{W_n}\int_{[1,\infty)}\int_0^1
\mathbb E^{x,u,t}[d_x^2]\,dt\,\mu(du)\,dx \le
|W_n|
\left(
 \Lambda_-\mathbb E U
 +\Lambda_-^2\mathbb E[U^2]
\right) \le C|W_n|.
\end{aligned}
\] 
To control the second factor in \eqref{44441}, choose \(S>2R\). There is a constant
\(C_{R,S}<\infty\) such that
\begin{equation}\label{eq:lambda_ineq}
\sup_{0\le r\le R}
\lambda^4e^{2r\lambda}
\le
C_{R,S}\big(e^{S\lambda}+e^{-S\lambda}\big)
\end{equation}
for every \(\lambda\in\mathbb R\). Indeed, for \(\lambda\ge0\), the
exponential \(e^{S\lambda}\) absorbs \(\lambda^4e^{2R\lambda}\) because
\(S>2R\). For \(\lambda<0\), we have \(e^{2r\lambda}\le1\), and
\(e^{-S\lambda}\) absorbs the polynomial factor.
Since \(A_n\) commutes with \(e^{rA_n}\), we have
\(
B_r^2=A_n^4e^{2rA_n}.
\)
Applying the inequality \eqref{eq:lambda_ineq} to the eigenvalues of \(A_n\)
gives
\[
\sup_{0\le r\le R}\Tr(B_r^2)
\le
C_{R,S}
\big(
\Tr(e^{SA_n})+\Tr(e^{-SA_n})
\big).
\]
Therefore, \eqref{eq:uniform-exponential-trace}, applied with \(S\) in place
of \(R\), implies
\(
\sup_{n\ge1}\sup_{0\le r\le R}
|W_n|^{-1}\mathbb E[\Tr(B_r^2)]
<\infty.
\)
Inserting this and 
\(
\sup_{n\ge1}
|W_n|^{-1}\mathbb E[\Tr(D_n^2)]
<\infty\)
 into \eqref{44441} 
yields
\eqref{eq:uniform-degree-trace}.

\end{proof}

We now prove the variance estimate for exponentially twisted linear
statistics in Proposition \ref{main prop}. The proof follows the  reduction of
\cite{hirsch-nam-otto-rgg}; the all-walk series estimate \cite[Lem. 5.3]{hirsch-nam-otto-rgg} is replaced by the
 trace bounds of
Lemma~\ref{lem:exponential_trace_bounds}. We sketch the proof and refer the reader to the proof of \cite[Thm. 5.2]{hirsch-nam-otto-rgg} for details.

\begin{proof}[Proof of Proposition \ref{main prop}]
Set
\(
    H_n(f)
    :=
    \Tr[f(A_n)e^{cA_n}] .
\)
For \(\check x\in\check W_n\),  
we identify \(A_n\) with the adjacency matrix of \(G_\emptyset(\xi_n,(\check x))\), and  we
write
\(
    \widehat A_n:=A_n\oplus0
\)
for the matrix corresponding to adding one isolated vertex at \(\check x\) to $G(\xi_n)$. Use the same ordering of the old vertices in
$A_n^{\check x}$ and $\widehat A_n$, with the inserted vertex
placed in the final matrix coordinate.
This matrix-coordinate convention is independent of its
reserved label in the edge-mark array.
In these coordinates, the upper-left block of
$A_n^{\check x}$ is exactly $A_n$.

Since \(f(0)=0\), the isolated vertex contributes
\(f(0)e^0=0\) to the trace, and therefore,
\(
    \Tr[f(\widehat A_n)e^{c\widehat A_n}]
    =
    \Tr[f(A_n)e^{cA_n}] .
\)
Thus,
\(
    \Delta_{\check x;(\check x)}H_n(f)
    =
    \Tr[f(A_n^{\check x})e^{cA_n^{\check x}}]
    -
    \Tr[f(\widehat A_n)e^{c\widehat A_n}] .
\)
Let
\(
    B_{\check x}
    :=
    A_n^{\check x}-\widehat A_n .
\)
By the construction of \(G_{\{1\}}(\xi_n,\check x)\), the induced graph on $\PP_n$ is exactly \(G_\emptyset(\xi_n,(\check x))\). Hence,
\(B_{\check x}\) is supported only on the row and column corresponding to the
inserted vertex \(\check x\). Let \(\mathbf e_{\check x}\) be the coordinate vector of the inserted vertex. Since $\widehat A_n\mathbf e_{\check x}=0$, we have
\(
    B_{\check x}\mathbf e_{\check x}
    =A_n^{\check x}\mathbf e_{\check x}
    =\mathbf a_{\check x}
\), where this vector is the same as the above defined $\mathbf a_{\check x}$. It is supported on the neighbors of \(\check x\) in \(G_{\{1\}}(\xi_n,\check x)\), and
\(
    B_{\check x}
    =
    \mathbf e_{\check x}\mathbf a_{\check x}^{\top}
    +
    \mathbf a_{\check x}\mathbf e_{\check x}^{\top}.
\)
For fixed \(n\), we have \(H_n(f)\in L^2\). Indeed, since
\(H^2(\mathbb R)\subseteq C_b(\mathbb R)\)  and
\(\|A_n\|_{\mathrm{op}}\le N_n\), we have
\(
    |H_n(f)|
    \le
    \|f\|_\infty N_n e^{|c|N_n},
\)
whose square has finite expectation because \(N_n\) is Poisson. Hence, the
Poincar\'e inequality from Theorem \ref{thm:poincare} gives
\(
    \Var(H_n(f))
    \le
    \int_{\check W_n}
    \mathbb E\big[(\Delta_{\check x;(\check x)}H_n(f))^2\big]\lambda(\d\check x).
\)

For $f\in L^1(\mathbb R)$, we use the Fourier transform
\(
\widehat f(\xi):=\frac{1}{2\pi}\int_{\mathbb R}f(x)e^{-i\xi x}\,\mathrm dx,
\)
so that $f(x)=\int_{\mathbb R}\widehat f(\xi)e^{i\xi x}\,\mathrm d\xi$ whenever Fourier inversion applies. To ensure this, assume first that \(f\in L^1(\mathbb R)\cap H^2(\mathbb R)\), and put
\[
    I_c(f)
    :=
    \int_{\mathbb R}
    |\widehat f(\xi)|\big(|\xi|+|c|\big)\,\d\xi .
\]
For \(s=c+i\xi\), Fourier inversion gives
\(
    \Tr[f(A_n^{\check x})e^{cA_n^{\check x}}]
    -
    \Tr[f(\widehat A_n)e^{c\widehat A_n}]
    =
    \int_{\mathbb R}
    \widehat f(\xi)
    \Tr\big[
    e^{sA_n^{\check x}}-e^{s\widehat A_n}
    \big]\,\d\xi .
\)
As in \cite[Eq. (50)]{hirsch-nam-otto-rgg}, Duhamel's formula yields  that
\(
    \Tr\big[
    e^{sA_n^{\check x}}-e^{s\widehat A_n}
    \big]
    =
    s\int_0^1
    \Tr\big[
    e^{tsA_n^{\check x}}
    B_{\check x}
    e^{(1-t)s\widehat A_n}
    \big]\,\d t .
\)
Using the representation of \(B_{\check x}\) and the rank-one trace identity $\Tr(Muv^{\top}N)=v^{\top}NMu$, we obtain
\(
\Tr[e^{tsA_n^{\check x}}B_{\check x}e^{(1-t)s\widehat A_n}]
=\mathbf a_{\check x}^{\top}e^{(1-t)s\widehat A_n}e^{tsA_n^{\check x}}\mathbf e_{\check x}
+\mathbf e_{\check x}^{\top}e^{(1-t)s\widehat A_n}e^{tsA_n^{\check x}}\mathbf a_{\check x}.
\)
Cauchy's inequality therefore gives
\[
    \big|
    \Tr[
    e^{tsA_n^{\check x}}
    B_{\check x}
    e^{(1-t)s\widehat A_n}
    ]
    \big|
    \le{}
    \|e^{tsA_n^{\check x}}\mathbf e_{\check x}\|_2
    \|\mathbf a_{\check x}^{\top}
    e^{(1-t)s\widehat A_n}\|_2
+
    \|e^{tsA_n^{\check x}}\mathbf a_{\check x}\|_2
    \|e^{(1-t)s\widehat A_n}\mathbf e_{\check x}\|_2 .
\]
Since \(A_n^{\check x}\) and \(\widehat A_n\) are real symmetric, the unitary
parts \(e^{it\xi A_n^{\check x}}\) and \(e^{i(1-t)\xi\widehat A_n}\) do not
change Euclidean norms. The vectors \(\mathbf e_{\check x}\) and
\(\mathbf a_{\check x}\) have non-negative entries, and both
adjacency matrices have non-negative entries. Hence, for either of these
vectors \(v\), the exponential series gives the entrywise bound
\(
    |(e^{tcA}v)_i|\le (e^{t|c|A}v)_i.
\)
Moreover, \(\widehat A_n\mathbf e_{\check x}=0\). Thus, the preceding display
is bounded by
\(
    X_t(\check x)^{1/2}\overline Y_{1-t}(\check x)^{1/2}
    +
    Y_t(\check x)^{1/2},
\)
where
\[
    X_t(\check x)
    :=
    \big\langle
    e^{2t|c|A_n^{\check x}}\mathbf e_{\check x},
    \mathbf e_{\check x}
    \big\rangle ,\qquad
    \overline Y_t(\check x)
    :=
    \big\langle
    e^{2t|c|\widehat A_n}\mathbf a_{\check x},
    \mathbf a_{\check x}
    \big\rangle , \qquad
    Y_t(\check x)
    :=
    \big\langle
    e^{2t|c|A_n^{\check x}}\mathbf a_{\check x},
    \mathbf a_{\check x}
    \big\rangle .
\]
Consequently, using \(|s|\le |\xi|+|c|\) together with Jensen's inequality in the
\(t\)-variable, we obtain the analogue of
\cite[Proof of Theorem~5.2, Step~1]{hirsch-nam-otto-rgg}:
\(
    \mathbb E\big[(\Delta_{\check x;(\check x)}H_n(f))^2\big]
    \le
    C I_c(f)^2
    \int_0^1
    \mathbb E\Big[
        X_t(\check x)\,\overline Y_{1-t}(\check x)
        +
        Y_t(\check x)
    \Big]\,\d t .
\)
Since \(A_n^{\check x}\) and \(\widehat A_n\) have non-negative entries and
\(\widehat A_n\le A_n^{\check x}\) entry-wise, the exponential series gives
\(e^{\alpha \widehat A_n}\le e^{\alpha A_n^{\check x}}\) entry-wise for every \(\alpha  \ge0\).
As \(\mathbf a_{\check x}\) has non-negative entries, it follows
that
\(
    \overline Y_t(\check x)
    \le
    Y_t(\check x)
\). We next bound the two matrix elements directly. Cauchy's inequality gives
\(
    X_t(\check x)^2
    \le
    \big\langle
        e^{4t|c|A_n^{\check x}}\mathbf e_{\check x},
        \mathbf e_{\check x}
    \big\rangle .
\)
Moreover, with \(\mathbf a_{\check x}=A_n^{\check x}\mathbf e_{\check x}\), we have
\(
    Y_t(\check x)^2
    \le
    \big\langle
        e^{4t|c|A_n^{\check x}}\mathbf a_{\check x},
        \mathbf a_{\check x}
    \big\rangle
    \left\langle
        \mathbf a_{\check x},\mathbf a_{\check x}
    \right\rangle .
\)
Since \(4t|c|\le4|c|\),
Lemma~\ref{lem:exponential_trace_bounds}, with \(R=4|c|\), yields 
\[
\sup_{0\le t\le1}\int_{\check W_n}\mathbb E[X_t(\check x)^2]\lambda(\d\check x)
+
\sup_{0\le t\le1}\int_{\check W_n}\mathbb E[Y_t(\check x)^2]\lambda(\d\check x)
\in O(|W_n|).
\]
Combining these two estimates and using Cauchy's inequality gives
\[
\begin{split}
    \int_{\check W_n}
    \mathbb E\big[
        X_t(\check x)\overline Y_{1-t}(\check x)
    \big]\lambda(\d\check x)
    &\le
    \int_{\check W_n}
    \mathbb E\big[
        X_t(\check x)Y_{1-t}(\check x)
    \big]\lambda(\d\check x)
    \le
    \sqrt{
    \int_{\check W_n}
    \mathbb E[X_t(\check x)^2]\lambda(\d\check x)}
    \sqrt{
    \int_{\check W_n}
    \mathbb E[Y_{1-t}(\check x)^2]\lambda(\d\check x)}
\end{split}
\]
which is $O(|W_n|)$.
Moreover, since \(\lambda=\d x\otimes\mu\) on \(\check W_n\) and \(\mu\) is
a probability measure,
\(
    \lambda(\check W_n)
    =
    |W_n|\mu([1,\infty))
    =
    |W_n|.
\)
Hence, Cauchy's inequality and the second-moment estimate for \(Y_t\) give
\[
    \int_{\check W_n}
    \mathbb E[Y_t(\check x)]\lambda(\d\check x)
    \le
    \sqrt{\lambda(\check W_n)}
    \sqrt{
    \int_{\check W_n}
    \mathbb E[Y_t(\check x)^2]\lambda(\d\check x)}
    \in O(|W_n|).
\]
Therefore,
\(
    \int_{\check W_n}
    \mathbb E\big[(\Delta_{\check x;(\check x)}H_n(f))^2\big]\lambda(\d\check x)
    \le
    C|W_n|I_c(f)^2 .
\)
By Cauchy--Schwarz and Plancherel's identity,
\[
\begin{aligned}
I_c(f)^2
&=
\left(
 \int_{\mathbb R}|\widehat f(\xi)|(|\xi|+|c|)\,d\xi
\right)^2
\le
\left(
 \int_{\mathbb R}(1+\xi^4)|\widehat f(\xi)|^2\,d\xi
\right)
\left(
 \int_{\mathbb R}
 \frac{(|\xi|+|c|)^2}{1+\xi^4}\,d\xi
\right) \le
C \bigl(\|f\|_2^2+\|f''\|_2^2\bigr).
\end{aligned}
\] 
Combining this estimate with the Poincar\'e inequality proves
the assertion for
$f\in L^1(\mathbb R)\cap H^2(\mathbb R)$ with $f(0)=0$.

Now let $f\in H^2(\mathbb R)$ satisfy $f(0)=0$.
Choose $\chi\in C_c^\infty(\mathbb R)$ with $0\le\chi\le1$ and
$\chi=1$ on $[-1,1]$, and set
\[
f_R(x):=\chi(x/R)f(x),\qquad R\ge1.
\]
Then $f_R\in L^1(\mathbb R)\cap H^2(\mathbb R)$ and $f_R(0)=0$.
Dominated convergence shows that $\|f_R-f\|_{H^2(\mathbb R)}\rightarrow0.$  
Sobolev embedding also gives $\|f_R-f\|_\infty\to0$.
For each fixed $n$, use $\|A_n\|_{\mathrm{op}}\le N_n$ to obtain
\[
\begin{aligned}
&\mathbb E\left|
 \operatorname{Tr}[f_R(A_n)e^{cA_n}]
 -
 \operatorname{Tr}[f(A_n)e^{cA_n}]
\right|^2 \le
\|f_R-f\|_\infty^2
\mathbb E\!\left[N_n^2e^{2|c|N_n}\right]
\longrightarrow0.
\end{aligned}
\]
The expectation on the right is finite as $N_n$ is Poisson.
Thus the trace statistics converge in $L^2$ for each fixed $n$,
and their variances converge.
Passing to $R\to \infty$ in the previous bound for $f_R$ yields
\(
\operatorname{Var}\bigl(
 \operatorname{Tr}[f(A_n)e^{cA_n}]
\bigr)
\le
C|W_n|\bigl(\|f\|_2^2+\|f''\|_2^2\bigr).
\)
The constant is the same as for $f_R$, and thus independent
of $n$. 
\end{proof}
\bigskip

\subsection*{Statement on AI use}

The original research ideas and the first drafts of all sections were our own. During subsequent revisions, we used AI, more precisely the latest versions of ChatGPT, to identify possible mathematical errors and suggest corrections. Moreover,  our original draft needed an exponential moment assumption that was stronger than the currently used \eqref{cond:EM} assumption. Here, AI suggested a preliminary version of Lemma 16 that allows us to prove the results in the stated generality. After that, we carefully went through this auxiliary result and revised its formulation so as to make it accessible.
All mathematical statements, proofs, computations, references, and conclusions were independently checked and verified by the authors, who take full responsibility for the content of the
manuscript.

\bibliographystyle{abbrv}
\bibliography{lit}

@book{last2018lectures,
  title={Lectures on the Poisson process},
  author={Last, G{\"u}nter and Penrose, Mathew},
  volume={7},
  year={2018},
  publisher={Cambridge University Press}
}

@article{DeprezWuthrich2019,
  author  = {Deprez, Philippe and W{\"u}thrich, Mario V.},
  title   = {Scale-Free Percolation in Continuum Space},
  journal = {Commun. Math. Stat.},
  volume  = {7},
  number  = {3},
  pages   = {269--308},
  year    = {2019},
  doi     = {10.1007/s40304-018-0142-0}
}

@article{BroidoClauset2019,
  author  = {Broido, Anna D. and Clauset, Aaron},
  title   = {Scale-free networks are rare},
  journal = {Nat. Commun.},
  volume  = {10},
  pages   = {1017},
  year    = {2019},
  doi     = {10.1038/s41467-019-08746-5}
}

@article{cltsparse,
  title={Central limit theorems for linear statistics of heavy tailed random matrices},
  author={Benaych-Georges, Florent and Guionnet, Alice and Male, Camille},
  journal={Comm. Math. Phys.},
  volume={329},
  number={2},
  pages={641--686},
  year={2014},
  publisher={Springer}
}

@article {mal_stab,
    AUTHOR = {Lachi\`eze-Rey, R. and Schulte, M. and Yukich, J. E.},
     TITLE = {Normal approximation for stabilizing functionals},
   JOURNAL = {Ann. Appl. Probab.},
    VOLUME = {29},
      YEAR = {2019},
    NUMBER = {2},
     PAGES = {931--993}
}

@article {mehler,
    AUTHOR = {Last, G. and Peccati, G. and Schulte, M.},    
     TITLE = {Normal approximation on {P}oisson spaces: {M}ehler's formula,
              second order {P}oincar\'{e} inequalities and stabilization},
   JOURNAL = {Probab. Theory Related Fields},
    VOLUME = {165},
      YEAR = {2016},  
    NUMBER = {3-4},
     PAGES = {667--723}
}

@article {shch1,
    AUTHOR = {Shcherbina, M. and Tirozzi, B.},
     TITLE = {Central limit theorem for fluctuations of linear eigenvalue
              statistics of large random graphs},
   JOURNAL = {J. Math. Phys.},
    VOLUME = {51},
      YEAR = {2010},
    NUMBER = {2},
     PAGES = {023523, 20}
}

@article{pcPerc,
	Author = {Penrose, M. D.},
	Journal = {Adv. Appl. Probab.},
	Number = {3},
	Pages = {536--556},
	Title = {On a continuum percolation model},
	Volume = {23},
	Year = {1991}
}

@article {nestmann,
    AUTHOR = {Last, G. and Nestmann, F. and Schulte, M.},
     TITLE = {The random connection model and functions of edge-marked {P}oisson processes: second order properties and normal
              approximation},
   JOURNAL = {Ann. Appl. Probab.},
    VOLUME = {31},
      YEAR = {2021},
    NUMBER = {1},
     PAGES = {128--168}
}

@article {trinh2,
    AUTHOR = {Can, V. H. and Trinh, K. D.},
     TITLE = {Random connection models in the thermodynamic regime: central
              limit theorems for add-one cost stabilizing functionals},
   JOURNAL = {Electron. J. Probab.},
    VOLUME = {27},
      YEAR = {2022},
     PAGES = {Paper No. 36, 40}
}

@article{shi-malik,
  author  = {Jianbo Shi and Jitendra Malik},
  title   = {Normalized Cuts and Image Segmentation},
  journal = {IEEE Trans. Pattern Anal. Mach. Intell.},
  volume  = {22},
  number  = {8},
  pages   = {888--905},
  year    = {2000},
  doi     = {10.1109/34.868688}
}

@inproceedings{ng-jordan-weiss,
  author    = {Andrew Y. Ng and Michael I. Jordan and Yair Weiss},
  title     = {On Spectral Clustering: Analysis and an Algorithm},
  booktitle = {Advances in Neural Information Processing Systems 14},
  pages     = {849--856},
  year      = {2002}
}

@article{von-luxburg,
  author  = {Ulrike von Luxburg},
  title   = {A Tutorial on Spectral Clustering},
  journal = {Stat. Comput.},
  volume  = {17},
  number  = {4},
  pages   = {395--416},
  year    = {2007},
  doi     = {10.1007/s11222-007-9033-z}
}

@book{penrose-rgg,
  author    = {Mathew Penrose},
  title     = {Random Geometric Graphs},
  series    = {Oxford Studies in Probability},
  volume    = {5},
  publisher = {Oxford University Press},
  year      = {2003}
}

@article{barthelemy-spatial,
  author  = {Marc Barth{\'e}lemy},
  title   = {Spatial Networks},
  journal = {Phys. Rep.},
  volume  = {499},
  number  = {1--3},
  pages   = {1--101},
  year    = {2011},
  doi     = {10.1016/j.physrep.2010.11.002}
}

@article{furedi-komlos,
  author  = {Zolt{\'a}n F{\"u}redi and J{\'a}nos Koml{\'o}s},
  title   = {The Eigenvalues of Random Symmetric Matrices},
  journal = {Combinatorica},
  volume  = {1},
  number  = {3},
  pages   = {233--241},
  year    = {1981},
  doi     = {10.1007/BF02579329}
}

@article{erdos-knowles-yau-yin,
  author  = {L{\'a}szl{\'o} Erd{\H{o}}s and Antti Knowles and Horng-Tzer Yau and Jun Yin},
  title   = {Spectral Statistics of {E}rd{\H{o}}s--{R}{\'e}nyi Graphs {I}: {L}ocal Semicircle Law},
  journal = {Ann. Probab.},
  volume  = {41},
  number  = {3B},
  pages   = {2279--2375},
  year    = {2013},
  doi     = {10.1214/11-AOP734}
}

@inproceedings{hamidouche-cottatellucci-avrachenkov,
  author    = {Mounia Hamidouche and Laura Cottatellucci and Konstantin Avrachenkov},
  title     = {Spectral Analysis of the Adjacency Matrix of Random Geometric Graphs},
  booktitle = {2019 57th Annual Allerton Conference on Communication, Control, and Computing},
  pages     = {208--214},
  year      = {2019},
  doi       = {10.1109/ALLERTON.2019.8919798}
}

@article{adhikari-adler-bobrowski-rosenthal,
  author  = {Adhikari, K. and  Adler, R. J. and Bobrowski, O. and Rosenthal, R.},
  title   = {On the Spectrum of Dense Random Geometric Graphs},
  journal = {Ann. Appl. Probab.},
  volume  = {32},
  number  = {3},
  pages   = {1734--1773},
  year    = {2022},
  doi     = {10.1214/21-AAP1720}
}

@article{hirsch-otto-svane,
  author  = {Christian Hirsch and Moritz Otto and Svane, A. M.},
  title   = {Normal Approximation for {G}ibbs Processes via Disagreement Couplings},
  journal = {Electron. J. Probab.},
  volume  = {30},
  year    = {2025},
  doi     = {10.1214/25-EJP1351}
}

@article{hirsch-nam-otto-rgg,
  author  = {Christian Hirsch and Kyeongsik Nam and Moritz Otto},
  title   = {Central Limit Theorem for Linear Eigenvalue Statistics of Random Geometric Graphs},
  journal = {arXiv preprint arXiv:2602.17006},
  year    = {2026}
}

@article{heerten-hirsch-otto,
  author  = {Heerten, N. and Hirsch, C. and Otto, M.},
  title   = {Moderate Deviations for Weight-Dependent Random Connection Models},
  journal = {J. Appl. Probab.},
  volume  = {62},
  number  = {4},
  pages   = {1406--1428},
  year    = {2025},
  doi     = {10.1017/jpr.2025.13}
}

@article{avrachenkov-bobu-dreveton,
  author  = {Avrachenkov, K. and Bobu, A. and Dreveton, M.},
  title   = {Higher-Order Spectral Clustering for Geometric Graphs},
  journal = {J. Fourier Anal. Appl.},
  volume  = {27},
  pages   = {22},
  year    = {2021},
  doi     = {10.1007/s00041-021-09825-2}
}

@incollection{bogdan,
  title={L{\'e}vy systems and moment formulas for mixed {P}oisson integrals},
  author={Bogdan, Krzysztof and Rosi{\'n}ski, Jan and Serafin, G. and Wojciechowski, {\L}.},
  booktitle={Stochastic Analysis and Related Topics: A Festschrift in Honor of Rodrigo Ba{\~n}uelos},
  pages={139--164},
  year={2017},
  publisher={Springer}
}

@book{Bhatia1997,
  author    = {Rajendra Bhatia},
  title     = {Matrix Analysis},
  series    = {Graduate Texts in Mathematics},
  volume    = {169},
  publisher = {Springer},
  address   = {New York},
  year      = {1997},
  doi       = {10.1007/978-1-4612-0653-8}
}

@article{pabst2026betti,
  author  = {Pabst, D.},
  title   = {Betti numbers in the Random Connection Model for higher-dimensional
             simplicial complexes and the {B}oolean model},
  journal = {Stochastic Process. Appl.},
  volume  = {201},
  pages   = {105071},
  year    = {2026},
  doi     = {10.1016/j.spa.2026.105071}
}

@article{gracarGrauerLuechtrathMoerters2019,
  author = {Gracar, Peter and Grauer, Arne and L{\"u}chtrath, Lukas and M{\"o}rters, Peter},
  title = {The age-dependent random connection model},
  journal = {Queueing Syst.},
  volume = {93},
  number = {3--4},
  pages = {309--331},
  year = {2019},
  doi = {10.1007/s11134-019-09625-y}
}

@article{gracarLuechtrathMoench2025,
  author = {Gracar, Peter and L{\"u}chtrath, Lukas and M{\"o}nch, Christian},
  title = {Finiteness of the percolation threshold for inhomogeneous long-range models in one dimension},
  journal = {Electron. J. Probab.},
  volume = {30},
  pages = {Paper No. 134, 29},
  year = {2025},
  doi = {10.1214/25-EJP1399}
}

@article {komjathy2,
    AUTHOR = {Gracar, P. and Heydenreich, M. and M\"{o}nch, C.
              and M\"{o}rters, P.},
     TITLE = {Recurrence versus transience for weight-dependent random
              connection models},
   JOURNAL = {Electron. J. Probab.},
    VOLUME = {27},
      YEAR = {2022},
     PAGES = {Paper No. 60, 31}
}

@book{netCom,
	Address = {Cambridge},
	Author = {Franceschetti, M. and Meester, R.},
	Publisher = {Cambridge University Press},
	Title = {Random Networks for Communication},
	Year = {2007}
}

\end{document}